\documentclass[12pt]{article}
\usepackage[utf8]{inputenc}

\title{\textbf{\Large Invariant foliations for endomorphims of $\mathbb{P}^2$ with a pluripotentialist product structure}}
\author{Virgile Tapiero}
\date{}

\usepackage[T1]{fontenc}
\usepackage[english, french]{babel}
\usepackage{hyperref}
\usepackage{amsmath, amssymb, amscd, amsthm}
\usepackage[all]{xy} 
\usepackage{esint}
\usepackage{color}

\newtheorem{thm}{Theorem}[section]
\newtheorem{lemme}[thm]{Lemma}
\newtheorem{prop}[thm]{Proposition}
\newtheorem{cor}[thm]{Corollary}
\newtheorem{defn}[thm]{Definition}
\newtheorem{rmq}[thm]{Remark}

\newcommand{\reels}{\mathbb{R}}
\newcommand{\cmplex}{\mathbb{C}}
\newcommand{\LLog}{\ \mathrm{Log}}

\providecommand{\keywords}[1]
{
  {\small \textbf{\textit{Keywords---}} #1}
}

\providecommand{\acknowledgements}[1]
{
  {\textbf{\textit{Acknowledgements---}} #1}
}

\newenvironment{changemargin}[2]{\begin{list}{}{%
\setlength{\topsep}{0pt}%
\setlength{\leftmargin}{0pt}%
\setlength{\rightmargin}{0pt}%
\setlength{\listparindent}{\parindent}%
\setlength{\itemindent}{\parindent}%
\setlength{\parsep}{0pt plus 1pt}%
\addtolength{\leftmargin}{#1}%
\addtolength{\rightmargin}{#2}%
}\item }{\end{list}}

\usepackage[left=3.20cm,right=3.20cm,top=2cm,bottom=2cm]{geometry}

\begin{document}
\hypersetup{pdfborder=0 0 0}

\maketitle

\begin{center}
    \begin{changemargin}{0.000005cm}{0.000005cm}
        \begin{otherlanguage}{english}
            \begin{abstract}
            {\footnotesize Let $f$ be a holomorphic endomorphism of $\mathbb{P}^2$, let $T$ be its Green current and $\mu=T\wedge T$ be its equilibrium measure. We prove that if $\mu$ has a local product structure with respect to $T$ then (an iterate of) $f$ preserves a local foliation $\mathcal{F}$ on a neighborhood of $\mathrm{Supp}(T)\backslash\mathcal{E}$, where $\mathcal{E}$ denotes the exceptional set of $f$. If the local foliation $\mathcal{F}$ extends through $\mathcal{E}$, then it extends to $\mathbb{P}^2$ and is an invariant pencil of lines. 
            } 
            \end{abstract}
        \end{otherlanguage}
    \end{changemargin}
\end{center}

\section{Introduction}

This article concerns holomorphic endomorphisms of $\mathbb{P}^2$ preserving geometric structures. Endomorphisms preserving pencils of curves, webs or foliations were classified by Dabija-Jonsson \cite{DabJon07,DabJon10} and Favre-Pereira \cite{FavPer09,FavPer15} by using algebraic geometry. Jonsson studied in \cite{jon99} the dynamics of polynomial skew products on $\cmplex^2$ by using pluripotential theory. He proved that the equilibrium measure $\mu$ satisfies a skew product property with respect to the Green current $T$ and provided formulas for the Lyapunov exponents. In particular, the two exponents of $\mu$ are larger than the logarithm of the algebraic degree $d$ of the endomorphism. Recently, Dupont-Taflin \cite{DT} extended Jonsson's results to endomorphisms of $\mathbb{P}^2$ preserving a pencil of lines.

In this article we are interested  in the reverse property. 
Let $f$ be an endomorphism of $\mathbb{P}^2$ of degree $d\geq2$, let $T$ be its Green current and $\mu:=T\wedge T$ be its  equilibrium measure.
We prove that if $\mu$ locally satisfies a product structure with respect to $T$, then there exists a holomorphic foliation near the support of $T$ which is $f$-invariant. Moreover, if the foliation extends to $\mathbb{P}^2$, then this extension is an invariant pencil of lines. The precise statements are given below.

Let us recall the definition of $T$ and $\mu$, we refer to Dinh and Sibony \cite{dinsib,sib99} for a general account. We have $T:=\lim_{n}d^{-n}{f^n}^*\omega_{\mathbb{P}^2}$, where $\omega_{\mathbb{P}^2}$ is the normalized Fubini-Study form, it is a positive closed
$(1,1)$-current on $\mathbb{P}^2$ satisfying $f^*T = d\times T$.  The product $\mu=T\wedge T$ is well defined since $T$ has local continuous potentials, it is a mixing invariant probability measure on $\mathbb{P}^2$ satisfying $f^*\mu=d^2\times \mu$. Its Lyapunov exponents satisfy $\lambda_1\geq\lambda_2\geq\frac{1}{2}\ \mathrm{Log}\ d$, as shown by Briend-Duval \cite{BriDuv99}.

For our purpose, let us develop the situation when $f$ preserves a pencil of lines as above.
It has the form $f[z:w:t]=[P(z,w):Q(z,w):R(z,w,t)]$ for a convenient choice of coordinates. Dupont-Taflin {\cite[Thm. 1.1]{DT}} proved that 
\begin{equation}\label{eq:Dupont-Taflin-formula}
    \mu = T\wedge \pi^*\mu_{\theta}\ \mathrm{and}\ \pi_*(\mu)=\mu_{\theta},
\end{equation}
where $\theta:=[P(z,w):Q(z,w)]$, $\pi : [z:w:t] \mapsto [z:w]$ (undefined at $[0:0:1]$) and $\mu_{\theta}$ is the equilibrium measure of $\theta$.
If moreover $\theta$ is a Lattès map of ${\mathbb{P}^1}$, then $$\mu=T\wedge\pi^*\mu_{\theta}\ll T\wedge \omega_{\mathbb{P}^2}$$ and the smallest exponent of $\mu$ is equal to $\lambda_{\theta} = \frac{1}{2}\ \mathrm{Log}\ d$, see {\cite[Cor. 1.3 and 1.4]{DT}}. Note that, without assuming the existence of an invariant pencil of lines, Dujardin proved that $\mu\ll T\wedge\omega_{\mathbb{P}^2}$ implies $\lambda_2=\frac{1}{2}\ \mathrm{Log}\ d$, see {\cite[Thm. 3.6]{Duj12}}  (see also {\cite{duptap23}} for another argument). 

We recall that Lattès maps of ${\mathbb{P}^2}$ are characterized by the three equivalent properties : $\lambda_1=\lambda_2=\frac{1}{2}\ \mathrm{Log}\ d$, $\mu\ll\omega_{\mathbb{P}^2}\wedge\omega_{\mathbb{P}^2}$, and $T$ is equal to a positive smooth $(1,1)$-form on a non empty open set. We refer here to the works by Berteloot, Dupont and Loeb \cite{berdup05, berloe98,BL01}. Similar (adapted) characterizations are also valid on $\mathbb{P}^1$. 

In the present article, we first begin by specifying the formula $\mu=T\wedge \pi^*\mu_{\theta}$ near most of repelling periodic points $a\in\mathrm{Supp}(\mu)$ when $\theta$ is Latt\`es: there exists Poincaré-Dulac coordinates $(Z_a,W_a)$ such that
\begin{equation}\label{eq:Dupont-Taflin-foormulaforrepellingcyclestheoremA}
    \mu = T\wedge dd^c|W_a|^2\ \mathrm{on}\ (\mathbb{P}^2,a) .
\end{equation}
This is done in Section \ref{sec:dupont-taflinformulaforrepulsivecycles}, by using results by Berteloot, Dupont, Loeb and Molino.

We are then motivated to study mappings satisfying Formula \eqref{eq:Dupont-Taflin-foormulaforrepellingcyclestheoremA} for at least one repelling periodic point. We introduce the Radon-Nikodym decomposition of the trace measure of $T$  with respect to $\mu$: $$\sigma_T :=T\wedge\omega_{\mathbb{P}^2} =\mu^a + \mu^s , $$ where $\mu^a\ll\mu$ and $\mu^s\perp \mu$.

\begin{thm}\label{thm:Main}
Let $f$ be a holomorphic map of $\mathbb{P}^2$ of degree $d\geq2$, $T$ be its Green current, $\mu=T\wedge T$ be its equilibrium measure and $\mathcal{E}$ be its exceptional set. Assume that $\mathrm{Supp}(\mu)\cap\mathcal{E}=\emptyset$. Assume moreover that 
\begin{enumerate}
    \item there exist a $N-$periodic repelling point $a\in\mathrm{Supp}(\mu)$ and Poincaré-Dulac coordinates $(Z_a,W_a)$ for $f^N$ at the point $a$ such that 
    $$\mu=T\wedge dd^c|W_a|^2\ \mathrm{ on }\ (\mathbb{P}^2,a)$$
    \item there exists an open set $\Omega$ charged by $\mu$ such that, for any open subset $V\subset\Omega$, $\mu(V)>0$ implies $\mu^s(V)>0$.
\end{enumerate}
Then there exists an open neighborhood $\mathcal{V}$ of $\mathrm{Supp}(T)\backslash\mathcal{E}$, and there exists a holomorphic foliation $\mathcal{F}$ defined on $\mathcal{V}$, such that $\mathcal{F}$ is invariant by $f^N$ on $\mathcal{V}$.
\end{thm}

Examples of mappings satisfying the hypothesis of this theorem are given by $f=[P(z,w):Q(z,w):t^d]$, where $\theta=[P(z,w):Q(z,w)]$ is Lattès. Indeed, in this case $\mathcal{E}$ is composed of the center of the pencil and the invariant line at infinity, we refer to Proposition \ref{prop:TheoremA} and \cite[§8]{duptap23} 
to check the items \textit{1.} and \textit{2.} More generally, if $f$ preserves a pencil of lines, then $\mathcal{E}$ contains its center. For a generic map of $\mathbb{P}^2$, $\mathcal{E}$ is empty.

The local first assumption $\mu = T\wedge dd^c|W_a|^2$ implies that $\mu\ll\sigma_T$ on $\mathbb{P}^2$, see Proposition \ref{prop:pullbackofduponttaflinformulabysigma}. In particular $f$ has a minimal Lyapunov exponent according to Dujardin’s theorem \cite{Duj12}. Note that the second assumption  implies that $f$ is not Lattès.

The proof of the Theorem \ref{thm:Main} requires several steps provided in Sections \ref{sec:linearizationofholomorphicgerms}, \ref{sec:dimensionofmeasureandapplication} and \ref{sec:proofoftheoremBconstructionofafoliation}. In Section \ref{sec:linearizationofholomorphicgerms} we introduce \textit{Poincaré maps} $\sigma: \cmplex^2 \to \mathbb{P}^2$ associated to Poincaré-Dulac coordinates $(Z_a,W_a)$. Roughly speaking, the leaves of the foliation $\mathcal{F}$ will be obtained by the images of the horizontal lines of $\cmplex^2$ by a Poincar\'e map $\sigma$.  To implement this idea, we prove in Section \ref{sec:dimensionofmeasureandapplication} a Patching Theorem : it allows, under the second assumption of the Theorem \ref{thm:Main}, to patch holomorphic $1-$forms $dW_1$ and $dW_2$ satisfying $T\wedge dd^c|W_i|^2\ll \mu$.

In Section \ref{sec:proofoftheoremBconstructionofafoliation}, we construct the foliation $\mathcal{F}$ near $\mathrm{Supp}(T)\setminus\mathcal{E}$, thus proving Theorem \ref{thm:Main}. First we explain that the local formula $\mu=T\wedge dd^c|W_a|^2$ near $a$ can be lifted to $\cmplex^2$ and gives a global formula $\sigma^*\mu=\sigma^*T\wedge dd^c|w|^2$ on $\cmplex^2$. Then in Section \ref{sec:finpreuvetheoremeB} we explain how to build the desired foliation $\mathcal{F}$ using Theorem \ref{thm:foliatedJ'}, proved in Section \ref{sec:proofoftheorem62}. The Theorem \ref{thm:foliatedJ'} asserts the following: we can construct a foliation $\mathcal{F}$, such that $\sigma^*\mathcal{F}$ coincide with the horizontal foliation of $\cmplex^2$, on a neighborhood of any compact set $J\subset \mathrm{Supp}(T)\backslash \mathcal{E}$ such that $J\supset \mathrm{Supp}(\mu)$ and 
$J=\mathrm{Supp}(\sigma_T|_J)$.\\

We obtain the following corollary.
\begin{cor}\label{cor:C} Let $f$ be an endomorphism of $\mathbb{P}^2$ satisfying the assumptions of $\mathrm{Theorem\ \ref{thm:Main}}$. Let $\mathcal{F}$ be a $f^N-$invariant holomorphic foliation on a neighborhood $\mathcal{V}$ of $\mathrm{Supp}(T)\backslash\mathcal{E}$ given by this theorem. 

Assume that for each point $p$ of $\mathrm{Supp}(T)\cap \mathcal{E}$ there exists an open neighborhood $U_p$ of $p$ such that $\mathcal{F}$ extends on $\mathcal{V}\cup U_p$. Then $\mathcal{F}$ uniquely extends to $\mathbb{P}^2$ and this extension is a pencil of lines invariant by $f^N$.
\end{cor}

Let us outline the proof, we give details in Section \ref{sec:recallonholomorphicfoliationsandsteinmanifolds}. Denote $\mathcal{U}:=\cup_{p\in \mathrm{Supp}(T)\cap \mathcal{E}} U_p$. We use the following extension theorem : any holomorphic foliation defined on an open neighborhood of a connected compact subset whose complementary set in $\mathbb{P}^2$ is Stein, extends to a holomorphic foliation on $\mathbb{P}^2$. This result is due to Lins Neto {\cite{AIF_1999__49_4_1369_0}}, see also Canales \cite{Can17}. Since the support of $T$ is a connected compact subset of $\mathbb{P}^2$ whose complementary set is Stein (see Fornaess-Sibony {\cite{forsib94}}, Ueda {\cite{10.2969/jmsj/04630545}}), we can apply Lins Neto’s theorem to the compact set $\mathrm{Supp}(T)\subset \mathcal{V}\cup\mathcal{U}$. We obtain that $\mathcal{F}$ extends to $\mathbb{P}^2$. We can then check that the extension is also invariant. Finally, we use Favre-Pereira classification \cite{FavPer09} of invariant foliations on $\mathbb{P}^2$.\\

\noindent\acknowledgements{The author thanks C. Dupont for numerous advice and suggestions concerning the writing of this article. He also thanks T.-C. Dinh for having suggested to add Lemma \ref{rmq:remarqueDinh2}, and for pointing out that the original Theorem \ref{lemme:patchinglemma} of the author could be generalized as in Remark \ref{sec:generalisationofthepatchingtheorem}. The author also thanks R. Dujardin for several fruitful discussions. This work was conducted within the France 2030 framework programme, Centre Henri Lebesgue ANR-11-LABX-0020-01.}

\section{Holomorphic foliations and invariance}\label{sec:recallonholomorphicfoliationsandsteinmanifolds}

\subsection{Holomorphic foliations on \texorpdfstring{$\mathbb{P}^2$}{TEXT}}\label{sec:holomorphicfoliations}

Let $U$ be an open subset of  $\mathbb{P}^2$. We define a foliation $\mathcal{F}$ on $U$ by a 
collection $(U_i,\omega_i)_{i\in I}$, where $\cup_{i\in U}U_i=U$ is an open cover, and where the $\omega_i\in\Omega^1(U_i)\backslash\{0\}$ are holomorphic $1-$forms, such that the following compatibility condition is satisfied:
\begin{equation}\label{eq:conditionintegrabiliteI}
    \omega_i\wedge \omega_j \equiv 0\ \mathrm{on}\ U_{i}\cap U_j.
\end{equation}
Two collections $(U_i,\alpha_i)_i$ and $(V_j,\beta_j)_j$ define the same foliation on $U$ if $\alpha_i\wedge \beta_j \equiv 0\ \mathrm{on}\ U_{i}\cap V_j$ whenever $U_i\cap V_j\neq\emptyset$. It is an equivalence relation. 

Let us assume that each $U_i$ is equipped with holomorphic coordinates $(Z_i,W_i)$ and that $\omega_i=f_idZ_i+g_idW_i$, where $f_i,g_i\in\mathcal{O}(U)$. Let us define the set $U^0:=U\backslash\left(\cup_{i\in I}\{f_i=0\}\cap\{g_i=0\}\right)$. Since we work in complex dimension $2$, the condition (\ref{eq:conditionintegrabiliteI}) ensures that the distribution of complex lines $p\in U^0\mapsto \mathrm{Ker}(\omega_i)_p$ is integrable in the sense of Frobenius, and thus there exists a sub-bundle $T\mathcal{F}\subset TU^0$ such that $T_p\mathcal{F}=\mathrm{Ker}(\omega_i)_p$, for any $p\in U^0$ and for any $i\in I$ such that $p\in U_i$.  The vector field $$v_i:=g_i\frac{\partial}{\partial Z_i} - f_i\frac{\partial}{\partial W_i}$$ satisfies 
$\omega_i(v_i) = 0\ \mathrm{on}\ U_i$ and so $T_p\mathcal{F} = \mathrm{Ker}(\omega_i)_p = \cmplex\cdot v_i(p)$ for any $p\in U_i\cap U^0$.

For the next Lemma, we introduce holomorphic $1$-forms $\eta_i=(h_idZ_i+k_idW_i)$ defining a foliation $\mathcal{G}$ on $U$. Let $(w_i)_i$ be the corresponding vector fields and $U^1:=U\backslash\left(\cup_{i\in I}\{h_i=0\}\cap\{k_i=0\}\right)$.

\begin{lemme}\label{lemma:F_1=F_2onD^2(1)thenF_1=F_2onD^2(2)}  For every $i\in I$, we have
\begin{enumerate}
    \item $\mathcal{F}|_{U_i}=\mathcal{G}|_{U_i}\ \Longleftrightarrow\ (\omega_i\wedge\eta_i = 0\ \mathrm{on}\ U_i)\ \Longleftrightarrow\  (\omega_i(w_i)=0\ \mathrm{on}\ U_i)$.
    \item  $\mathcal{F}|_{U_i}=\mathcal{G}|_{U_i}\ \Longleftrightarrow\ \left(T_p\mathcal{F}=T_p\mathcal{G},\ \forall p\in U_i\cap U^0\cap U^1\right)$.
    \item Assume $U=U_i,\ \forall i\in I,$ and that $U$ is connected. If there exists a non empty open subset $V\subset U$ such that $\mathcal{F}|_V=\mathcal{G}|_V$, then by analytic continuation $\mathcal{F}=\mathcal{G}$.
\end{enumerate}
\end{lemme}
\noindent\underline{\textbf{Proof :}} \ \\
\noindent\textit{1.} The first equivalence is our definition of foliations. For the second equivalence, just observe that $\omega_i\wedge\eta_i = \left|^{f_i}_{g_i}\phantom{0}^{h_i}_{k_i}\right|dZ_i\wedge dW_i$ and that $\omega_i(w_i)=f_ik_i-g_ih_i=\left|^{f_i}_{g_i}\phantom{0}^{h_i}_{k_i}\right|$.\\

\noindent\textit{2.} We have to prove the reverse implication. Let us assume that $T_p\mathcal{F}=T_p\mathcal{G}$ for any $p\in U_i\cap U^0\cap U^1$. Since $T_p\mathcal{F}=\mathrm{Ker}(\omega_i)_p$ and $T_p\mathcal{G}=\cmplex\cdot w_i(p)$, we get $\omega_i(w_i)(p)=0$ for any $p\in U_i\cap U^0\cap U^1$. Since $U_i$ is connected, $\omega_i(w_i)(p)=0$ for every $p\in U_i$ by analytic continuation. Hence $\mathcal{F}|_{U_i}=\mathcal{G}|_{U_i}$ according to the first item.\\

\noindent\textit{3.} If $\mathcal{F}|_V=\mathcal{G}|_V$ on $V\subset U$, then $\omega_i\wedge\eta_i = 0$ on $V$. Since $\omega_i\wedge\eta_i=\left|^{f_i}_{g_i}\phantom{0}^{h_i}_{k_i}\right|dZ_i\wedge dW_i$, we deduce that $\left|^{f_i}_{g_i}\phantom{0}^{h_i}_{k_i}\right|=0$ on $V$ and thus on $U$ by analytic continuation. It means that $\omega_i\wedge\eta_i=0$ on $U_i$ for every $i\in I$, and thus $\mathcal{F}=\mathcal{G}$ on $U$.\qed

\subsection{Pull-back of foliations and invariance}\label{sec:invariantfoliations}

Let $U$ and $V$ be two open subsets of $\mathbb{P}^2$ and let $f:V\to U$ be a holomorphic surjective map which is not constant on each connected component of $V$. Let $\mathcal{F}$ be a foliation on $U$ defined by a collection $(U_i,\omega_i)_{i\in I}$. For each $i\in I$, the $1-$form $f^*\omega_i$ is not null on each connected component of $f^{-1}(U_i)$ and $f^{*}\omega_i\wedge f^{*}\omega_j= f^{*}(\omega_i\wedge \omega_j)=0$ on $f^{-1}(U_i)\cap f^{-1}(U_j)$. Then the collection $(U_i,f^*\omega_i)_{i\in I}$ defines a foliation $f^*\mathcal{F}$ on $V$, called the pull-back of $\mathcal{F}$ by $f$.

\begin{defn}\label{defn:invariancedunfeuilletagelocal} Assume $V=U$ and that $f^{-1}(V)\cap V\neq\emptyset$. We say that $\mathcal{F}$ is invariant by $f$ on $V$ if $f^*\mathcal{F}=\mathcal{F}$ on ${f^{-1}(V)\cap V}$. 
\end{defn}
For $V=U=\mathbb{P}^2$ the pairs $(\mathcal{F},f)$ where $f^*\mathcal{F}=\mathcal{F}$, have been classified by Favre-Pereira. We shall need the following  version of their classification.

\begin{thm}[Favre-Pereira {\cite{FavPer09}}]\label{thm:Favre-Pereira} Let $f:\mathbb{P}^2\rightarrow\mathbb{P}^2$ be a holomorphic map of degree $d\geq2$. Let $\mathcal{F}$ be a foliation on $\mathbb{P}^2$ invariant by $f$. Then in appropriate homogeneous coordinates $[x:y:z]$ on $\mathbb{P}^2$, one of the following cases holds:
\begin{enumerate}
    \item $\mathcal{F}$ is the pencil on lines  $\pi[z:w:t]=[z:w]$. In this case $f=[P(x,y):Q(x,y):R(x,y,z)]$ with $P,Q,R$ homogeneous polynomials of degree $d$. 
    \item $\mathcal{F}$ is not a pencil of lines through a point in $\mathbb{P}^2$ and $f$ has the form:
    \begin{enumerate}
        \item[i.] $f=[x^d:y^d:z^d]$ or $f=[z^d:x^d:y^d]$ or
        \item[ii.] $f=[x^d:y^d:R(x,y,z)]$ or $f=[y^d:x^d:S(x,y,z)]$, where $R$ and $S$ depends on the three variables $(x,y,z)$.
    \end{enumerate}
\end{enumerate}
We refer to {\cite{FavPer09}} for a description of the foliations $\mathcal{F}$, and for the polynomials $R,\ S$.
\end{thm}

In particular, Latt\`es maps on $\mathbb{P}^2$ do not preserve any foliation, however they can preserve webs \cite{FavPer15}. We shall use Theorem \ref{thm:Favre-Pereira} in the following form.

\begin{cor}\label{cor:fpreservesFonP2thenF=pencil} If $f$ preserves a foliation $\mathcal{F}$ on $\mathbb{P}^2$ and if the smallest Lyapunov exponent of $\mu$ satisfies $\lambda_2=\frac{1}{2}\ \mathrm{Log}\ d$, then in appropriate homogeneous coordinates $\mathcal{F}$ is a pencil of lines given by the fibers of $\pi[z:w:t]=[z:w]$.
\end{cor}

We then can prove the following result.

\begin{prop}\label{prop:extensionisalsoinvariant}
Let $f$ be a holomorphic map of $\mathbb{P}^2$ of degree $d\geq2$ and let $\mathcal{F}$ be a $f-$invariant foliation defined on a neighborhood of a fixed point of $f$. If $\mathcal{F}$ admits an extension $\mathcal{F}'$ on $\mathbb{P}^2$, then $\mathcal{F}'$ is invariant by $f$. Moreover, if  $\lambda_2=\frac{1}{2}\ \mathrm{Log}\ d$ then, in appropriate coordinates, $\mathcal{F}'$ is the pencil of lines given by $\pi[z:w:t]=[z:w]$. In this case the rational map induced by $f$ on $\mathbb{P}^1$ is Lattès.
\end{prop}
\noindent\textbf{\underline{Proof :}} The second part of the statement is a direct consequence of Corollary \ref{cor:fpreservesFonP2thenF=pencil}. 
Let $a\in\mathbb{P}^2$ be the fixed point given by the statement. $\mathcal{F}'$ coincide with $\mathcal{F}$ on a neighborhood of $a$, and $f^*\mathcal{F}=\mathcal{F}$ near $a$. So $f^*\mathcal{F}'=f^*\mathcal{F}=\mathcal{F}$ near $a$. Using again that $\mathcal{F}'$ is an extension of $\mathcal{F}$, we deduce that $f^*\mathcal{F}'=\mathcal{F}'$ near $a$ and thus on $\mathbb{P}^2$ by analytic continuation (Lemma \ref{lemma:F_1=F_2onD^2(1)thenF_1=F_2onD^2(2)}).\qed 

\section{Poincaré-Dulac coordinates}\label{sec:dupont-taflinformulaforrepulsivecycles}

Let us consider $f:(\mathbb{P}^2,a)\to(\mathbb{P}^2,a)$ a holomorphic germ such that the eigenvalues $(\chi_1,\chi_2)$ of $d_af$ satisfy $|\chi_1|\geq |\chi_2|>1$. 
By Poincaré-Dulac theorem, there exists a germ of biholomorphism $\sigma_0:(\cmplex^2,0)\to(\mathbb{P}^2,a)$ and a polynomial mapping $D(z,w)=(\chi_1z+cw^q,\chi_2w)$, with $q\geq2$, such that the following commutative relation holds: 
$$f\circ\sigma_0=\sigma_0\circ D\ \mathrm{on}\ (\cmplex^2,0).$$
The map $\sigma_0^{-1}$ induces holomorphic coordinates $\sigma_0^{-1}=(Z_a,W_a)$, we call them \textit{Poincaré-Dulac coordinates} for the germ $f$ of $(\mathbb{P}^2,a)$.

\subsection{Latt\`es mappings on \texorpdfstring{$\mathbb{P}^1$}{TEXT} and multipliers on \texorpdfstring{$\mathbb{P}^2$}{TEXT}}\label{sec:Berteloot-LoebtheoremonLattèsmap}

We begin by recalling some results due to Berteloot-Loeb. For a rational map $\theta$ of degree $d\geq2$  on $\mathbb{P}^1$, we use the following definitions: 
\begin{enumerate}
    \item[\textbullet] A \textit{regular point} $p\in\mathbb{P}^1$ is a point such that $\mu_{\theta}$ is a strictly positive smooth $(1,1)-$form on an open neighborhood of $p$ in $\mathbb{P}^1$. Let $\Omega_{\theta}$ be the set of regular points and let $A_{\theta}:=\mathbb{P}^1\backslash\Omega_{\theta}$.
    \item[\textbullet] $R_1(\theta)$ is the set of repelling fixed points of $\theta$, it is contained in $\mathrm{Supp}(\mu_{\theta})$.
\end{enumerate}
We refer to {\cite[Theorem 1.1 \& Proposition 2.5]{berloe98}} and {\cite[Proposition 4.1]{BL01}}.

\begin{thm}[Berteloot-Loeb]\label{thm:BLI}\ 
\begin{enumerate}
    \item If $\theta$ is a Lattès map then $A_{\theta}$ is a finite set, otherwise $A_{\theta} = \mathbb{P}^1$.
    \item Let $a\in R_1(\theta)$, $D(w) := \theta'(a) w$ and suppose that $a\in\Omega_{\theta}$. There exists an invertible holomorphic germ $\sigma_0:(\cmplex,0)\rightarrow(\mathbb{P}^1,a)$ such that the following relation holds:
\begin{equation}\label{eq:diagrammedecommutationdeBLThmI}
    \sigma_0^{-1}=:W_0\ \mathrm{and}\ W_0\circ \theta = D\circ W_0\ \mathrm{on}\ (\cmplex,0)
\end{equation}
and such that $\mu_{\theta}=dd^c|W_0|^2$ (which implies $|\theta'(a)|=\sqrt{d}$). We call the coordinate $W_0$ a Poincaré-Dulac coordinate for $\theta$ at the point $a$.
\end{enumerate}
\end{thm}

Now we deal with the multipliers of repelling cycles for mappings $f$ of $\mathbb P^2$ of degree $d$, of equilibrium measure $\mu$ and of Lyapunov exponents $\lambda_1\geq\lambda_2$. Let $\varepsilon > 0$ and $R_{\mu}^{n,\varepsilon}$ be the set of $n$-periodic repelling points $a\in\mathrm{Supp}(\mu)$ satisfying 
$$\frac{1}{n}\LLog\left|\mathrm{det}_{\cmplex}\ d_af^n\right|\geq \lambda_1+\lambda_2 - 2\varepsilon . $$

\begin{prop}[Berteloot-Dupont-Molino {\cite[Lemma 4.5]{bdm07}}]\label{bdmm}\ 
\begin{enumerate}
\item Let $V\subset\mathbb{P}^2$ be an open set satisfying $\mu(V)>0$. There exists $n_{\varepsilon}\geq1$ such that for any $n\geq n_{\varepsilon}$:
\begin{equation}\label{eq:estimationBriendDuvallocalisee}
    \mathrm{Card}(R_{\mu}^{n,\varepsilon}\cap V)\geq d^{2n}(1-\varepsilon)^3\mu(V).
\end{equation}
In particular  $R_{\mu}^{\varepsilon}:=\bigcup_{n\geq1}R^{n,\varepsilon}_{\mu}$ is dense in $\mathrm{Supp}(\mu)$.
\item Let $R^{n}_{\mu}$ denote the set of $n$-periodic repelling points satisfying 
$$\frac{1}{n}\LLog\left|\mathrm{det}_{\cmplex}\ d_af^n\right|>\LLog\ d.$$
If $\lambda_1>\lambda_2$ then  $R_{\mu}:=\bigcup_{n\geq1}R^{n}_{\mu}$ is  dense in $\mathrm{Supp}(\mu)$.
\end{enumerate}
\end{prop}

We note that (\ref{eq:estimationBriendDuvallocalisee}) is actually proved in {\cite[Lemma 4.5]{bdm07}} for $V=\mathbb{P}^2$, but the same proof works  for a general open set $V$ charged by $\mu$. The second item is a direct consequence of the first item. 

\subsection{Endomorphisms of \texorpdfstring{$\mathbb{P}^2$}{TEXT} preserving a pencil of lines}{\label{sec:endomorphismspreservingapenciloflinesandproofofproposition23}}

We recall some facts on mappings preserving a pencil of lines through a point in $\mathbb{P}^2$. 
Let $f$ be a holomorphic map of $\mathbb{P}^2$ of degree $d\geq2$, and let $T$ and $\mu$ being the Green current of $f$ and $\mu$ its equilibrium measure. We denote $\lambda_1\geq\lambda_2$ the Lyapunov exponents of $\mu$. Let us assume that $f$ preserves the pencil of lines given $\pi:[z:w:t]\mapsto[z:w]$. We denote $\theta$ the rational map such that $\pi\circ f = \theta\circ \pi$ on $\mathbb{P}^2\backslash\{[0:0:1]\}$.

Recall that $A_{\theta}$ is the complementary set of regular points $p$ ($\mu_{\theta}$ is smooth on a neighborhood of $p$) in $\mathbb{P}^1$. We denote $E_{\theta}:= \pi^{-1}(A_{\theta})$, and we define $R_{\mu,\theta}$ the set of periodic and repelling points $a\in\mathrm{Supp}(\mu)\backslash E_{\theta}$ such that (with $n$ the period of $a$):  
$$\frac{1}{n}\LLog\left|\mathrm{det}_{\cmplex}\ d_af^n\right|>\LLog\ d.$$

\begin{prop}\label{prop:TheoremA}
    Assume that $\theta$ is Lattès.
    Then ${E_{\theta}}\cup\{[0:0:1]\}$ is a finite union of projectives lines which do not contain $\mathrm{Supp}(\mu)$, and $R_{\mu,\theta}$ is dense in $\mathrm{Supp}(\mu)\backslash E_{\theta}$. For any $a\in R_{\mu,\theta}$ of period $n\geq1$, there exists $(Z_a,W_a)$ Poincaré-Dulac coordinates for $f^n$ near $a$ such that :
    \begin{equation*}
        \mu = T\wedge dd^c|W_a|^2\ \mathrm{on}\ (\mathbb{P}^2,a).
    \end{equation*}
\end{prop}

Let us give the arguments for the proof. Note that the relation $\pi_*\mu = \mu_{\theta}$ (provided by Dupont-Taflin \eqref{eq:Dupont-Taflin-formula}) implies:
\begin{equation}\label{eq:pi(supp(mu))=supp(mu_theta)}
    \pi(\mathrm{Supp}(\mu))= \mathrm{Supp}(\mu_{\theta}) .
\end{equation}
According to Theorem \ref{thm:BLI}, the assumption $\theta$ Lattès implies that $E_{\theta}\cup\{[0:0:1]\}$ is a union of a finite number of projective lines passing through the point $[0:0:1]$. We have $[0:0:1]\not\in\mathrm{Supp}(\mu)$ (this point belongs to the Fatou set as a super attractive fixed point) and we can not have $\mathrm{Supp}(\mu)\subset E_{\theta}$ since $\pi(\mathrm{Supp}(\mu))=\mathrm{Supp}(\mu_{\theta})=\mathbb{P}^1$ (by \eqref{eq:pi(supp(mu))=supp(mu_theta)}) and since $\pi(E_{\theta})=A_{\theta}$ is finite. We deduce that $\mathrm{Supp}(\mu)\backslash E_{\theta}$ is a non empty open subset of $\mathrm{Supp}(\mu)$, and therefore 
Proposition \ref{bdmm} implies the first assertion of Proposition \ref{prop:TheoremA}. It remains to prove  the second one.

Let us fix $a\in\mathrm{Supp}(\mu)$ a repelling $N$-periodic point satisfying $\frac{1}{N}\LLog|\chi_1\chi_2|>\LLog\ d$, where $\chi_1 , \chi_2$ are the eigenvalues of $d_af^N$ (with $|\chi_1|\geq|\chi_2|$). Let $a_0:=\pi(a)$, since $a$ is a $N-$periodic point of $f$, the relation $\pi\circ f=\theta\circ\pi$ ensures that $a_0$ is a $N-$periodic point of $\theta$. This relation also ensures that $(\theta^N)'(a_0)$ is an eigenvalue of $d_af^N$. Finally, $a_0$ is a repelling fixed point of $\theta^N$.

We recall that $\Omega_{\theta}=\mathbb{P}^1\backslash A_{\theta}$. Now if we assume that $a\in R_{\mu,\theta}$ $i.e.$ $a_0\in\Omega_{\theta}$, then to prove the second assertion in Proposition \ref{prop:TheoremA}, it remains to prove the following:

\begin{lemme} With the preceding notations, if $a_0\in\Omega_{\theta}$ then 
\begin{enumerate}
\item $|\chi_1|>|\chi_2|=\sqrt{d^N},$
\item There exists Poincaré-Dulac coordinates $(Z_a,W_a)$ for $f^N$ such that :
$$\mu=T\wedge dd^c|W_a|^2\ \textrm{and}\ W_a\circ f^N=\chi_2 W_a \textrm{ on } (\mathbb{P}^2,a) . $$
\end{enumerate}
\end{lemme}
\noindent\textbf{\underline{Proof :}} Since $f$ and $f^N$ have same equilibrium measure $\mu$ and Green current $T$ (similarly for $\theta$ and $\theta^N$), we can assume that $N=1$. 
By assumption the fixed point $a_0$ belongs to $\Omega_{\theta}$, hence by Theorem \ref{thm:BLI} there exists an open neighborhood $U_0$ of $a_0$ and there exists a Poincaré-Dulac coordinate $W_0$ on $U_0$ such that $\mu_{\theta}=dd^c|W_0|^2$ on $U_0$. We denote $W_a:=W_0\circ \pi$ on $\pi^{-1}(U_0)\ni a$. Then injecting the normal form $\mu_{\theta}=dd^c|W_0|^2$ into Dupont-Taflin formula \eqref{eq:Dupont-Taflin-formula}, one gets $\mu = T\wedge dd^c|W_a|^2$ on $\pi^{-1}(U_0)$. 

At this stage, by the commutative relation \eqref{eq:diagrammedecommutationdeBLThmI}, we have:
\begin{equation}\label{eq:thmAstep3W_a=lambdaW_a}
    W_a\circ f = W_0 \circ \pi \circ f = W_{0}\circ \theta\circ\pi = \lambda\times W_0\circ\pi = \lambda W_a,
\end{equation}
where $\lambda=\theta'(a_0)$. According to Theorem \ref{thm:BLI}, it satisfies $|\lambda|=\sqrt{d}$.

Let $\widetilde{Z}_a$ be a submersion on an open neighborhood $U_a$ of $a$ such that $U_a\subset \pi^{-1}(U_0)$  and such that $\widetilde{\xi}_a:=(\widetilde{Z}_a,W_a)$ define holomorphic coordinates centered at $a$. If $\widetilde{D}_f:=\widetilde{\xi}_a\circ f\circ\ (\widetilde{\xi}_a)^{-1}$ and $p_w:(z,w)\mapsto w$, then 
$$p_w\circ \widetilde{D}_f=p_w\circ\widetilde{\xi}_a\circ f\circ (\widetilde{\xi}_a)^{-1} = W_a\circ f\circ (\widetilde{\xi}_a)^{-1} = \lambda W_a\circ (\widetilde{\xi}_a)^{-1}=\lambda p_w,$$
where the third equality is given by (\ref{eq:thmAstep3W_a=lambdaW_a}). This shows that $\widetilde{D}_f$ has the form
\begin{equation}\label{eq:thmAstep4TildeD_f=(h(z,w),lambdaw)}
    \widetilde{D}_f(z,w)=(h(z,w),\lambda w),
\end{equation}
where $h\in\mathcal{O}(\cmplex^2,0)$. By construction $d_0\widetilde{D}_f$ has the same eigenvalues $\chi_1$ and $\chi_2$ than $d_af$ which by assumption satisfy $|\chi_1\chi_2|>{d}$. From (\ref{eq:thmAstep4TildeD_f=(h(z,w),lambdaw)}) we deduce that $\lambda$ is also an eigenvalue of $d_0\widetilde{D}_f$: $\lambda\in\{\chi_1,\chi_2\}$. But  $|\lambda|=\sqrt{d}$, hence $\chi_2 =\lambda$ and $|\chi_1|>|\chi_2|$. By Lemma \ref{lemme:complitionintoPooincareDulaccoordinates} (proved below) there exist a chart $\xi_a$ such that $D_f\circ \xi_a = \xi_a\circ f$ on $(\cmplex^2,0)$, where $D_f(z,w)=(\chi_1z+cw^q,\lambda w)$ and $p_w\circ\xi_a=W_a$. Finally $\xi_a =(Z_a,W_a)$ are Poincaré-Dulac coordinates for $f$ at $a$ and $\mu = T\wedge dd^c|W_a|^2$ near $a$.\qed\\

\begin{lemme}\label{lemme:complitionintoPooincareDulaccoordinates} Let $g$ be a germ fixing $0$ in $\cmplex^2$ such that there exists $\xi=(Z,W)$ an invertible germ fixing $0$ such that $D:=\xi\circ g\circ \xi^{-1}$ has the form $D(z,w)=(h(z,w),\chi_2 w)$. We denote $ \chi_1 := \partial_z h(0,0)$ and we assume 
\begin{equation}\label{eq:conditionforthelemmaofsemilinearization}
    |\chi_1|>|\chi_2|>1.
\end{equation}
Then there exists $\widetilde{\xi} = (\widetilde{Z},W)$ an invertible germ fixing $0$ and there exists a germ $\widetilde{D}(z,w)=(\chi_1z+\widetilde{c}w^{\widetilde{q}},\chi_2 w)$ such that $\widetilde{\xi}\circ g\circ \widetilde{\xi}^{-1} = \widetilde{D}$.
\end{lemme}

\noindent\textbf{\underline{Proof :}} 
Since $|\chi_1|>|\chi_2|>1$, Poincaré-Dulac theorem yields the existence of $\xi'$ an invertible germ fixing $0$ such that:
\begin{equation}\label{eq:doublediagramme}
    D'\circ \xi' = \xi'\circ D\ \mathrm{on}\ (\cmplex^2,0)
\end{equation}
where
$$D'(z,w)=(\chi_1z+cw^q,\chi_2w)\ \mathrm{and}\ c=0\ \mathrm{if}\ \chi_2^q\neq \chi_1.$$
Let $W'$ be the second component of $\xi'$, our goal  is to show that $W'$ is given by $W':(z,w)\mapsto \gamma\times w$ on $(\cmplex^2,0)$ for some constant $\gamma\in\cmplex^*$. Taking $\widetilde{\xi}:=\xi''\circ\xi'\circ\xi$ with $\xi''(z,w) :=(z,\frac{1}{\gamma}w)$ and $\widetilde{D}:=\xi''\circ D' \circ (\xi'')^{-1}$, the desired conclusion follows.\\

As a holomorphic function, $W'$ is a power series near $0$ of the form :
$$W'(z,w) = \sum_{r\geq1}\left(\gamma_{r,0}z^r+\gamma_{r-1,1}z^{r-1}w+\cdots+\gamma_{0,r}w^r\right)=\sum_{r\geq1}\mathcal{O}'(r).$$
The function $h$ also expends into a power series near $0$ of the form :
$$h(z,w) = \sum_{r\geq1}\left(H_{r,0}z^r+H_{r-1,1}z^{r-1}w+\cdots+H_{0,r}w^r\right)= [H_{1,0}z + H_{0,1}w] + \mathcal{O}_{\geq 2}.$$
Of course $H_{1,0}$ is equal to $\chi_1$ (we recall that $\chi_1$, $\chi_2$ are the eigenvalues of $d_0 D$). The equality \eqref{eq:doublediagramme} provides 
\begin{equation}\label{eq:chi_2W'=W'(h,chi_2w)}
    \chi_2W'(z,w) = W'(h(z,w),\chi_2w)\ \mathrm{on }\ (\cmplex^2,0).
\end{equation}
Observe that (\ref{eq:chi_2W'=W'(h,chi_2w)}) implies (looking at terms of order $1$ in the power series):
$$\chi_2\gamma_{1,0}z + \chi_2\gamma_{0,1}w = (\gamma_{1,0}H_{1,0}z +\gamma_{1,0}H_{0,1}w) + \gamma_{0,1}\chi_2 w.$$
Since $H_{1,0}=\chi_1$ we deduce, looking at coefficients of the variable $z$, that $\chi_2\gamma_{1,0}=\gamma_{1,0}\chi_1$.  Since $\chi_1\neq\chi_2$, we get 
\begin{equation}\label{eq:stepr=1}
    \gamma_{1,0}=0\ \ i.e.\ \ \mathcal{O}'(1)=\gamma_{0,1}w .
\end{equation}

We can now prove by induction on $r\geq2$ that:
\begin{equation*}
    \forall k\in\{2,\cdots,r\},\ \mathcal{O}'(k)=0\ \ i.e.\ \ W'(z,w) = \gamma_{0,1}w + \sum_{k\geq r+1}\mathcal{O}'(k)\ \mathrm{on}\ (\cmplex^2,0).
\end{equation*}
Let us explain the base case $r=2$, the induction step $r>2$ can be treated analogously. Using (\ref{eq:chi_2W'=W'(h,chi_2w)}) and (\ref{eq:stepr=1}), we get :
$$\chi_2\gamma_{2,0}z^2 + \chi_2\gamma_{1,1}zw + \chi_2\gamma_{0,2}w^2 = \gamma_{2,0}(H_{1,0}z + H_{0,1}w)^2 + \gamma_{1,1}(H_{1,0}z + H_{0,1}w)(\chi_2w) + \gamma_{0,2}(\chi_2w)^2.$$
Since $\chi_1=H_{1,0}$ we deduce, looking at the coefficients of $z^2$, that $\chi_2\gamma_{2,0} = \gamma_{2,0}\chi_1^2$, it yields $\gamma_{2,0}=0$ by (\ref{eq:conditionforthelemmaofsemilinearization}). Using again $\chi_1=H_{1,0}$ and looking at the coefficient of $zw$ we obtain $\chi_2\gamma_{1,1}=\gamma_{1,1}\chi_1\chi_2$, thus $\gamma_{1,1}=0$ by (\ref{eq:conditionforthelemmaofsemilinearization}). At last looking at the coefficients of $w^2$ we obtain $\chi_2\gamma_{0,2}=\gamma_{0,2}\chi_2^2$, and thus $\gamma_{0,2}=0$. We have proved $\mathcal{O}'(2)=0$.\qed 

\section{Poincaré maps for repelling fixed points}\label{sec:linearizationofholomorphicgerms}
 
Let us fix a dilating germ  $f:(\mathbb{P}^2,a)\to(\mathbb{P}^2,a)$ and let us consider a germ $D(z,w) = (\chi_1z+cw^q,\chi_2w)$, with $c\in\cmplex$, $q \geq 2$, and Poincaré-Dulac coordinates $\sigma_0^{-1}=(Z_a,W_a)$ such that $f\circ\sigma_0=\sigma_0\circ D$ near $0$ on $\cmplex^2$. We assume that $f$ is induced by a holomorphic map of $\mathbb{P}^2$ of degree $d\geq 2$. 

From the local commutative relation $f\circ\sigma_0=\sigma_0\circ D$ near $0$ in $\cmplex^2$, we shall construct a global commutative relation
semi-conjugating $f$ to $D$ via a globally defined holomorphic map $\sigma:\cmplex^2\rightarrow\mathbb{P}^2$. The construction of $\sigma$ is classical, see for instance  Berteloot-Loeb \cite{BL01} when $D$ is linear.

Let $U_a$ be a connected open neighborhood of $a$ on which $\sigma_0^{-1}$ is defined. We denote $U_0:=\sigma_0^{-1}(U_a)$. Let $\mathbb{D}_{\varepsilon}^2 \subset U_0 \cap D^{-1}(U_0)$ be a small bidisc  such that $D^{-1}  (\mathbb{D}_{\varepsilon}^2) \subset \mathbb{D}_{\varepsilon}^2$. In particular the following diagram commutes  
\begin{equation}\label{eq:diagrammedelapropositiondeladeveloppante}
    \xymatrix{
    \mathbb{D}_{\varepsilon}^2 \ar@{->}[d]_{\sigma_0} \ar[r]^{D} & D(\mathbb{D}_{\varepsilon}^2) \ar@{->}[d]^{\sigma_0} \\
    \sigma_0(\mathbb{D}_{\varepsilon}^2) \ar@{->}[r]^{}^{f} & f (\sigma_0(\mathbb{D}_{\varepsilon}^2))}
\end{equation}

\begin{prop}[Poincar\'e maps]\label{prop:thedevolepingmap}  \
\begin{enumerate}
    \item There exists an increasing sequence of integers $(n_k)_k$ such that $D^{-n_k} (\mathbb{D}_k^2) \subset \mathbb{D}_{\varepsilon}^2$. In particular the restriction $\sigma_k$ of $f^{n_k} \circ\sigma_0\circ D^{-n_k}$ to $\mathbb{D}_k^2$ is well defined. It satisfies $f \circ \sigma_k = \sigma_k \circ D$ on $\mathbb{D}_k^2 \cap D^{-1} (\mathbb{D}_k^2)$.
    \item $\sigma_k$ is an open mapping and  $\# \sigma_k ^{-1}(p) \leq d^{2 n_k}$ for every $p\in\sigma_k(\mathbb{D}_k^2)$. 
    \item The critical values $E_k:=\sigma_k(\mathrm{Crit}\ \sigma_k)$ are included in $f^{n_k}(\mathrm{Crit}\ f^{n_k})$, which is an algebraic subset of codimension $1$ of $\mathbb{P}^2$.     
    \item For every positive integers $k\leq l$, $\sigma_k=\sigma_l$
    on $\mathbb{D}_k^2$. Hence one can define the holomorphic map $\sigma:\mathbb{C}^2\rightarrow\mathbb{P}^2$  by $\sigma :=\lim_{k}\sigma_k$.
    \item The relation $f\circ \sigma=\sigma\circ D$ holds on $\cmplex^2$.
    \item The map $\sigma$ is open and have discrete fibers on $\cmplex^2$. 
    \item If $E:=\sigma(\mathrm{Crit}\ \sigma)$ then $E=\bigcup_{k\geq0}E_k$ and $\mu(E)=\sigma_T(E)=0$. 
\end{enumerate}
\end{prop}
\noindent\textbf{\underline{Proof :}} We explain the first and the fourth items, the others being a direct consequence of the previous ones, or could be deduced from the fact that $f$ is a degree $d^2$ ramifying covering map of $\mathbb{P}^2$, and from the fact that $\mu$ and $\sigma_T$ do not charge proper analytic subsets of $\mathbb{P}^2$ (cf. \cite{dinsib} for details).\\

\noindent\textit{1.} We have $D^n = (\chi_1^n z + nc \chi_1^{n-1} w^q , \chi_2^n w)$ and $D^{-n} = (\chi_1^{-n} z - nc \chi_1^{-(n-1)} w^q , \chi_2^{-n} w)$. Hence, for every $k \geq 1$, there exists $n_k$ such that $D^{-n_k} (\mathbb{D}_k^2)
\subset \mathbb{D}_{\varepsilon}^2$. Let $\sigma_k$ be the restriction of $f^{n_k} \circ\sigma_0\circ D^{-n_k}$ to $\mathbb{D}_k^2$. One can  assume that $n_k$ is increasing. We have $f \circ \sigma_k = \sigma_k \circ D$ on $\mathbb{D}_k^2 \cap D^{-1} (\mathbb{D}_k^2)$ thanks to the commutative diagram (\ref{eq:diagrammedelapropositiondeladeveloppante}).\\

\noindent\textit{4.}
Let $y \in \cmplex^2$ and assume that $D^m(y)\in\mathbb{D}_{\varepsilon}^2$ for some $m \geq 1$. From the inclusion $D^{-1}  (\mathbb{D}_{\varepsilon}^2) \subset \mathbb{D}_{\varepsilon}^2$, we get $\{ y,D(y),\cdots,D^m(y) \} \subset \mathbb{D}_{\varepsilon}^2$.
Now if $p\in D^n(\mathbb{D}_{\varepsilon}^2)\cap D^{n+m}(\mathbb{D}_{\varepsilon}^2)$ with $p=D^n(x)=D^{n+m}(y)$, then $D^m(y)=x\in \mathbb{D}_{\varepsilon}^2$ and thus, as explained above, $D^j(y)\in\mathbb{D}_{\varepsilon}^2$ for any $j\in\{0,\cdots,m\}$. By the commutative diagram (\ref{eq:diagrammedelapropositiondeladeveloppante}) we get 
$$f^m\circ\sigma_0(y) = f^{m-1}\circ\sigma_0\circ D(y)=\cdots=\sigma_0(D^m(y))=\sigma_0(x).$$ 
Composing by $f^n$ we obtain $f^{n+m}\circ\sigma_0\circ D^{-(n+m)}(D^{n+m}(y)) = f^n\circ\sigma_0\circ D^{-n}(D^n(x))$. If $k\leq l$ then taking $n=n_k$ and $m=n_l-n_k\geq0$, one has $\sigma_{l}=\sigma_k$ on $D^{n_k}(\mathbb{D}_{\varepsilon}^2)\cap D^{n_l}(\mathbb{D}_{\varepsilon}^2)\supset \mathbb{D}_k^2$. It allows to define $\sigma:\cmplex^2\rightarrow\mathbb{P}^2$ by $\sigma|_{\mathbb{D}_k^2}:=\sigma_k$ for any $k\geq0$. It also satisfies $\sigma=\lim_k\sigma_k$.\qed\\

It is important to know whether $\sigma(\cmplex^2)$ contains $\mathrm{Supp}(\mu)$ or not. Let us introduce the exceptional set $\mathcal{E}$ of $f$: 
\begin{equation}\label{eq:exceptionalsetE(f)off}
    \mathcal{E}:=\mathbb{P}^2\backslash \left\{x\in\mathbb{P}^2 \ \left| \ \mu_{x,n} := \frac{1}{d^{2n}}\sum_{y\in f^{-n}(x)}\delta_y\underset{n\to+\infty}{\overset{\mathrm{weak-}\star}{\longrightarrow}}\mu \right.\right\}.
\end{equation}
The set $\mathcal{E}$ is the largest totally $f-$invariant proper algebraic subset of $\mathbb{P}^2$ \cite{dinsib}.

\begin{prop}\label{prop:imageofdevelopmentmapssigmacontainsP^2minusE(f)} If the fixed point $a$ belongs to the support of $\mu$, then 
$$ \mathrm{Supp}(\mu)\backslash\mathcal{E} \subset \mathbb{P}^2\backslash\mathcal{E} \subset 
 \sigma(\cmplex^2) . $$
\end{prop}
\noindent\textbf{\underline{Proof :}} Let $q\in\mathbb{P}^2\backslash\mathcal{E}$ and let $\rho\in C^0(\mathbb{P}^2,\reels^+)$ be a continuous function such that $\rho\equiv1$ on a ball $B_a$ centered at $a$ and such that $\mathrm{Supp}(\rho)\subset 2B_{a}$, we can assume that $\sigma_0^{-1}$ is defined on $2B_a$. By definition of $\mathcal{E}$ one has $\lim_n \langle\mu_{q,n},\rho\rangle=\langle\mu,\rho\rangle\geq \mu(B_a)>0$ because $a$ belongs to the support of $\mu$. Thus there exists $n_k\geq1$ such that $\langle\mu_{q,n_k},\rho\rangle>0$ and then, by definition of $\mu_{q,n_k}$, there exists a point $q_{n_k}\in f^{-n_k}(q)$ such that $\rho(q_{n_k})>0$. Since the support of $\rho$ is included in $2B_a$ we have $q_{n_k} \in 2B_a\subset\mathrm{Dom}(\sigma_0^{-1})$. This allows us to consider  $p_{n_k}:=D^{n_k}\circ\sigma_0^{-1}(q_{n_k})$, it satisfies $\sigma(p_{n_k})=f^{n_k}(q_{n_k})=q$.\qed\\

We conclude this section with the following lemma. This type of statement is classical and arises from Briend-Duval \cite{BriDuv99} work. We refer also to {\cite[Lemme 1.1.32]{dup02}} where the arguments can be borrowed to prove the lemma. 

\begin{lemme}\label{lemma:inversionlocale}
    Let $U$ be a bounded open subset of $\cmplex^2\backslash(\mathrm{Crit}\ \sigma)$. There exists $r>0$ satisfying the following property. Let $x\in U$ and let $p:=\sigma(x)$, then $B(p,r)\subset\sigma(U)$ and the
    inverse branche $\sigma_{x,p}^{-1}$ of $\sigma:U\to\sigma(U)$ sending $p$ to $x$ is defined on $B(p,r)$.
\end{lemme}

\begin{cor}\label{cor:relevementdecheminparsigma}
    Let $U$ as in $\mathrm{Lemma}$ $\ref{lemma:inversionlocale}$. For every Lipschitz path $\gamma:[0,1]\to \sigma(U)$ and for every $x_0\in U$ such that $\sigma(x_0)=\gamma(0)$, then there exists a (unique) continuous path $\Tilde{\gamma}:[0,1]\to U$ such that $\sigma\circ\Tilde{\gamma}=\gamma$. 
\end{cor}
\noindent\textbf{\underline{Proof:}} Lemma \ref{lemma:inversionlocale} ensures that there exists $r>0$ such that any inverse branches $\sigma_{x,p}^{-1}$ of $\sigma:U\to \sigma(U)$ (sending $p$ to $x$) is defined on ${B}(p,r)$. Let $N>\kappa/r$ be an integer, where $\kappa$ is the Lipschitz constant of $\gamma$. Put $p_j:=\gamma(j/N)$, $0\leq j\leq N-1$. Observe that $\gamma([j/N,(j+1)/N])\subset B(p_j,r)$ since $\gamma$ is $\kappa-$Lipschitz. Using inverse branches, we can defined a sequence of points $x_0,x_1,\cdots,x_{N-1}$ such that $\sigma(x_j)=p_j$. We can then define the desire path $\widetilde{\gamma}:t\in[0,1]\mapsto \gamma\circ\sigma_{x_j,p_j}^{-1}(t)$, $t\in[j/N,(j+1)/N]$.\qed

{\section{The Patching Theorem}\label{sec:dimensionofmeasureandapplication}}

We recall the Radon-Nikodym decomposition of $\sigma_T=T\wedge\omega_{\mathbb{P}^2}$ with respect to $\mu$ :
\begin{equation}\label{eq:decompositionofRNofsigma_Twithrespecttomu}
    \sigma_T=\mu^a + \mu^s.
\end{equation}
Here $\mu^a$ is absolutely continuous with respect to $\mu$ and $\mu^s\perp\mu$. Our aim is to prove:

\begin{thm}[Patching Theorem]\label{lemme:patchinglemma}
Let $U$ be a connected chart of $\mathbb{P}^2$ with two systems of holomorphic coordinates $(Z_1,W_1)$ and $(Z_2,W_2)$. Let us suppose that these systems of coordinates satisfy:
\begin{equation}\label{eq:conditionsTwedgeddc|W|2<<mu}
    T\wedge dd^c|W_i|^2\ll \mu\ \mathrm{on}\ U,\ i\in\{1,2\}.
\end{equation}
If moreover $\mu^s(U)>0$, then there exists $\beta\in\mathcal{O}^*(U)$ such that $dW_2=\beta dW_1\ \mathrm{on}\ U.$
\end{thm}

We give two proofs of the theorem. The first one exploits the hermitian properties of the Green current $T$ which can be seen as a positive singular (not smooth) metric on $\mathbb{P}^2$. More precisely, we express $T\wedge dd^c|W_2|^2$ as a sum of measures absolutely continuous with respect to the measures $T\wedge dd^c|Z_1|^2$ and $T\wedge dd^c|W_1|^2$. By using some results about dimension of measures, recalled in the section just below, we conclude that $dW_2$ and $dW_1$ are proportional. Variations of this technique can be used to detect fine properties of the current $T$, for example in \cite{duptap23} it is used to give a new proof of Dujardin's theorem: $\mu\ll\sigma_T$ implies $\lambda_2=\frac{1}{2}\ \mathrm{Log}\ d$. 

The second proof do not use specific properties of the Green current, except that it has local continuous potentials. In particular, we mention in Remark \ref{sec:generalisationofthepatchingtheorem} that this proof can be adapted for a more general setting. 

\subsection{Hausdorff dimension of measures}

Let $\nu$ be a finite Borel measure on $\mathbb{P}^2$. The \textit{pointwise lower dimension} of $\nu$ at $x\in\mathbb{P}^2$ is defined by
$$\underline{d}_{\nu}(x):=\liminf_{r\rightarrow0^+}\frac{\LLog\ \nu B(x,r)}{\LLog\ r}.$$
Since $T$ has continuous Hölder potentials, there exists a lower estimate for the lower pointwise dimension of $\sigma_T$ as follows.

\begin{prop}[Dinh-Sibony {\cite[Proposition 1.18]{dinsib}}]\label{prop:HolderpotentialsfortheGreenCurrent} For every $x \in \mathbb{P}^2$ we have $T=dd^c(u_x)$ on a small ball $B_x$, where $u_x$ is a plurisubharmonic function on $B_x$ which is $\gamma-$Hölder for any $\gamma\in]0,\gamma_0[$, where 
$$\gamma_0=\min\left\{1,\frac{\LLog(d)}{\LLog(d_{\infty})}\right\}>0\ \ \mathrm{and}\ \ d_{\infty}=\lim_{n\to+\infty} \left(\sup_{p\in \mathbb{P}^2}||d_pf^n||\right)^{1/n}.$$
In particular, $\underline{d}_{\sigma_T}(x)\geq 2+\gamma_0$ for every  $x\in\mathbb{P}^2$.
\end{prop}

\begin{cor}\label{eq:corollairevanishinglemma} 
Let $U$ be a connected open set of $\mathbb{P}^2$ such that $\nu(U)>0$. \ 
\begin{enumerate}
    \item Assume that there exists $\gamma>2$ such that $\underline{d}_{\nu}(x)\geq\gamma$ for $\nu-$almost every $x\in\mathbb{P}^2$. Then for any $h\in\mathcal{O}(U)$, $\nu\{x\in U:h(x)=0\}>0$ implies $h\equiv0$ on $U$. 
    \item Assume that $\nu\leq C\sigma_T$ for some constant $C>0$. Let $h\in\mathcal{O}(U)$ such that $h(p)=0$ for $\nu-$almost every point $p\in U$. Then $h\equiv0$ on $U$.
\end{enumerate}
\end{cor}
\noindent\textbf{\underline{Proof:}} \ \\
\noindent\textit{1.} 
Let $\mathrm{HD}(A)$ denote the Hausdorff dimension of $A \subset \mathbb{P}^2$. Young {\cite[Proposition 2.1]{you82}} proved that if $\underline{d}_{\nu}(x)\geq \gamma$ for $\nu-$almost every $x\in\mathbb{P}^2$, then $\mathrm{HD}(A)\geq\gamma$ for every Borel set $A$ of positive $\nu$-measure. Using the assumptions, we get  
\begin{equation}\label{eq:dim_H(Z(h))>2}
    \mathrm{HD}\left(\{x\in U:h(x)=0\}\right)\geq  \gamma>2.
\end{equation}
If $h$ were not null on $U$, then the analytic set $\{x\in U:h(x)=0\}$ would have Hausdorff dimension $2$, see for instance Chirka's book {\cite[Corollary 1 p.23]{chirka1989complex}}. But this is not compatible with (\ref{eq:dim_H(Z(h))>2}).

\noindent\textit{2.} 
Let $\widetilde{\nu}:=C\sigma_T-\nu\geq0$ which is a finite Borel measure, so we have for any $x\in \mathbb{P}^2$, $\underline{d}_{(\nu+\widetilde{\nu})}(x) = \underline{d}_{\sigma_T}(x)\geq 2+\gamma_0$ by Proposition \ref{prop:HolderpotentialsfortheGreenCurrent}. It follows that $\underline{d}_{\nu}(x)\geq\underline{d}_{\sigma_T}(x)$ for any $x\in\mathbb{P}^2$. By applying the first item with $\gamma:=2+\gamma_0$, we get the result.\qed

\subsection{Proof of the Patching Theorem}\label{sec:thepatchingtheorem}

\textbf{\textit{Proof using measure theory.---}} For any differentiable function $h$ on $U$, we use the following notations :
\begin{equation}\label{eq:definitionofpartilZ1partialW1}
    \frac{\partial h}{\partial Z_1}:=\frac{\partial}{\partial z}\left[h\circ \xi_1^{-1}\right]\circ \xi_1\ \mathrm{and}\ \frac{\partial h}{\partial W_1}:=\frac{\partial}{\partial w}\left[h\circ \xi_1^{-1}\right]\circ \xi_1,
\end{equation}
where $\xi_1=(Z_1,W_1)$ is the chart associated to the coordinates $(Z_1 ,W_1)$. We are going to prove that $dW_2=\frac{\partial W_2}{\partial W_1}dW_1$, the function $\frac{\partial W_2}{\partial W_1}$ being  defined by (\ref{eq:definitionofpartilZ1partialW1}) with $W_2 = p_w \circ \xi_2$. It is a local problem, it is sufficient to prove the equality on each $U'$ open with $\overline{U'}\subset U$. So we can assume without loss of generality that the holomorphic functions $Z_1, W_1, Z_2$ and $W_2$ are holomorphic on a neighborhood of $\overline{U}$ and thus the partial derivatives ${\partial W_2}/{\partial Z_1}$ and ${\partial W_2}/{\partial W_1}$ are bounded on $U$.

The Green current $T$ on $U$ can be written in the coordinates $(Z_1,W_1)$ as a $(1,1)-$differential form with complex measures coefficients :
$$T|_U = \sigma_1\frac{i}{2}dZ_1\wedge d\overline{Z}_1 + \left(\Lambda\frac{i}{2}dZ_1\wedge d\overline{W}_1\right) + \overline{\left(\Lambda\frac{i}{2}dZ_1\wedge d\overline{W}_1\right)} + \sigma_2\frac{i}{2}dW_1\wedge d\overline{W}_1,$$
with $\sigma_1$ and $\sigma_2$ the positive measures given by $\sigma_1 = T\wedge \frac{i}{2}dW_1\wedge d\overline{W}_1$ and $\sigma_2 = T\wedge \frac{i}{2}dZ_1\wedge d\overline{Z}_1$. The complex measure $\Lambda$ is given by
$\Lambda = - T\wedge\frac{i}{2}dW_1\wedge d\overline{Z}_1$. By Cauchy-Schwarz inequality, we have for every Borel set $A$ in $U$:
\begin{equation}\label{eq:Cauchy-Schwarzinequalityforpatchinglemma}
    |\Lambda(A)|\leq \sqrt{\sigma_1(A)}\sqrt{\sigma_2(A)}.
\end{equation}
The trace $\sigma_U:=\sigma_1 + \sigma_2$
gives a positive measure which is equivalent to the trace measure $\sigma_T$ on $U$:
\begin{equation}\label{eq:sigma_Tequivtomu+lambda}
    \exists C>0\ :\ \frac{1}{C}\sigma_T|_U\leq \sigma_U\leq C\sigma_T|_U
\end{equation}

Let us denote $\lambda:=\sigma_2\ \mathrm{and}\ \mu_1^a:=\sigma_1$. By (\ref{eq:conditionsTwedgeddc|W|2<<mu}) we have $\mu_1^a=\psi_1\mu|_U\ll \mu$, with $\psi_1\in L^1(\mu|_U)$. Let $\lambda = h_1\mu|_U + \mu^s_1$ 
be the Radon-Nikodym decomposition of $\lambda$ with respect to $\mu|_U$, where $h_1\in L^1(\mu|_U)$ is a non negative function, and $\mu^s_1$ is a positive measure on $U$ singular with respect to $\mu$. We have (recall that $\mu^s(U)>0$):
\begin{equation}\label{eq:equivalencedesmesuressinguliereslocaleetglobale}
    \frac{1}{C}\mu^s|_U\leq \mu^s_1\leq C\mu^s|_U\ \Longrightarrow \mu_1^s(U)>0
\end{equation}
with the same $C>0$ than in (\ref{eq:sigma_Tequivtomu+lambda}). Indeed, by using \eqref{eq:sigma_Tequivtomu+lambda} we have $\mu_1^s\leq C(\mu^a+\mu^s)$ and $\mu^s|_U\leq C[(h_1+\psi_1)\mu|_U+\mu^s_1]$, thus $\mu_1^s\leq C\mu^s$ and $\mu^s|_U\leq C\mu_1^s$ since $\mu_1^s\perp \mu^a$ and $\mu^s\perp\mu$. 

The idea now is to decompose  $dW_2$ in terms of the $1-$forms $dZ_1$ and $dW_1$:
$$dW_2 = \frac{\partial W_2}{\partial Z_1}dZ_1 + \frac{\partial W_2}{\partial W_1}dW_1.$$
Denoting $\alpha := \frac{\partial W_2}{\partial Z_1},\ \beta := \frac{\partial W_2}{\partial W_1}\ \mathrm{and}\ \gamma := \alpha\overline{\beta}$, we have :
$$dd^c|W_2|^2 = |\alpha|^2dd^c|Z_1|^2 + 2\mathrm{Re}\left[\gamma \frac{i}{2}dZ_1\wedge d\overline{W_1}\right] + |\beta|^2dd^c|W_1|^2.$$
Wedging by $T$ this equality, we get :
$$\underbrace{T\wedge dd^c|W_2|^2}_{=:\mu_2^a} = |\alpha|^2\underbrace{(T\wedge dd^c|Z_1|^2)}_{=\lambda} + \underbrace{2\mathrm{Re}\left[\gamma \left(T\wedge\frac{i}{2}dZ_1\wedge d\overline{W}_1\right)\right]}_{=:\Lambda'} + |\beta|^2\underbrace{(T\wedge dd^c|W_1|^2)}_{=\mu_1^a},$$
and thus we have:
\begin{equation}\label{eq:decompositiondeW_2enfontiondeZ_1etW_1}
    \mu_2^a = |\alpha|^2\lambda + \Lambda' + |\beta|^2\mu_1^a,
\end{equation}
with $\Lambda'$ a signed measure which satisfies for any Borel set $A\subset U$:
\begin{equation}\label{eq:majorationofLambda'byLambda}
    |\Lambda'(A)|\leq 2\left(\max_{U}|\gamma|\right)|\Lambda(A)|\leq 2\left(\max_{U}|\gamma|\right)\sqrt{\mu_1^a(A)}\sqrt{\lambda(A)},
\end{equation}
where the second inequality comes from (\ref{eq:Cauchy-Schwarzinequalityforpatchinglemma}). Let $A\subset U$ be any Borel set, then using (\ref{eq:decompositiondeW_2enfontiondeZ_1etW_1}) and using the fact that $\lambda = h_1\mu|_U + \mu_1^s$ we get :
\begin{align*}
    0\leq \int_A|\alpha|^2\ \mathrm{d}\mu^s_1 & = \mu_2^a(A) -  \int_A|\beta|^2\ \mathrm{d}\mu_1^a - \Lambda'(A) - \int_Ah_1|\alpha|^2\ \mathrm{d}\mu\\
                                        & \leq \mu_2^a(A) + |\Lambda'(A)|
                                        \leq \mu_2^a(A) + 2\left(\max_{U}|\gamma|\right)\sqrt{\mu_1^a(A)}\sqrt{\lambda(A)}\ \mathrm{by\ \textnormal{(\ref{eq:majorationofLambda'byLambda}})}.
\end{align*}
Since $\mu_1^a$ and $\mu_2^a$ are absolutely continuous with respect to $\mu$ and since $\mu_1^s\perp \mu$ on $U$, there exists $\mathcal{A}\subset U$ a Borel set of measure $0$ for $\mu_1^a$ and $\mu_2^a$ such that $\mu_1^s=\mu_1^s(\cdot\cap\mathcal{A})$. So using the preceding inequalities we have for any $A\subset U$:
\begin{align*}
    0\leq \int_A|\alpha|^2\ \mathrm{d}\mu^s_1 & = \int_{A\cap \mathcal{A}}|\alpha|^2\ \mathrm{d}\mu^s_1 \leq 0 + 2\left(\max_{U}|\gamma|\right)\sqrt{0}\sqrt{\lambda(A\cap\mathcal{A})} = 0,
\end{align*}
and thus $\int_A|\alpha|^2\ \mathrm{d}\mu_1^s=0$ for any $A\subset U$. The measure $\mu_1^s$ charges the open set $U$ by (\ref{eq:equivalencedesmesuressinguliereslocaleetglobale}), hence $\alpha(p)=0$ for $\mu_1^s-$almost every point $p\in U$. By construction $\mu_1^s\leq C\sigma_T$ for some $C>0$, so we can apply the second item of Corollary \ref{eq:corollairevanishinglemma} with $\nu:=\mu_1^s$ to conclude that $\alpha\equiv0$ on $U$. Thus we obtain $dW_2=\beta dW_1$ on $U$, and finally $\beta\in\mathcal{O}^*(U)$ since $W_1$ and $W_2$ are submersions on $U$.\qed\\

\textbf{\textit{Proof using pluripotential theory.---}} If we assume that $dW_1\wedge dW_2$ is not null on $U$, then $dW_1\wedge dW_2\neq0$ on a Zariski open set $V$ of $U$. In $V$, the coordinates $(W_1,W_2)$ create a chart and thus we have:
$$\mu^s\leq\sigma_T\ll T\wedge dd^c|W_1|^2 + T\wedge dd^c|W_2|^2\ll\mu\ \mathrm{by\ assumption}.$$
By the Chern-Levine-Nirenberg inequality, since $T$ has local continuous potentials, the current $T$ has no mass on $U\backslash V$ which is an analytic subset. Therefore the property $\mu^s\ll\mu$ on $V$ extends on $U$. Hence $\mu^s(U)=0$ and the result follows.\qed
\begin{rmq}\label{sec:generalisationofthepatchingtheorem}
    This second proof does not use the properties of the Green current $T$. The Patching Theorem thus can be generalized as follows.
    
    Let $X$ be a complex surface (a two dimensional complex manifold) equipped with a smooth $>0$ closed $(1,1)-$form $\omega$. Let $S$ be a positive closed $(1,1)-$current on $X$ with continuous local potentials and being of finite mass. Denote $\nu:=S\wedge S\ ,\ \sigma_S:=S\wedge\omega$.
    These wedge products are well defined since $S$ has continuous local potentials. Then $\nu$ is a positive finite Borel measure on $X$ (possibly identically null) and so is $\sigma_S$, so we can write the Radon-Nikodym decomposition of $\sigma_S$ with respect to $\nu$: $\sigma_S=\nu^a + \nu^s$, where $\nu^a\ll\nu$ and $\nu^s\perp \nu$. 

    The same arguments than above allows to prove the following. Let $U\subset X$ be a connected open subset of $X$ equipped with two systems of holomorphic coordinates $(Z_1,W_1)$ and $(Z_2,W_2)$. Assume that $S\wedge dd^c|W_i|^2\ll \nu$ for $i\in\{1,2\}$. If $\nu^s(U)>0$, then the holomorphic $1-$forms $dW_1$ and $dW_2$ patch together: there exists $\beta\in\mathcal{O}^*(U)$ such that $dW_2=\beta dW_1$ on $U$.
\end{rmq}

\section{Construction of a foliation near \texorpdfstring{$\mathrm{Supp}(T)\backslash\mathcal{E}$}{TEXT}}\label{sec:proofoftheoremBconstructionofafoliation}

We want to prove the Theorem \ref{thm:Main}. Let $a\in \mathrm{Supp}(\mu)$ be a repelling $N$-periodic point  of $f$ with Poincaré-Dulac coordinates $\sigma_0^{-1}=(Z_a,W_a)$. Let $D(z,w)=(\chi_1z+cw^q,\chi_2w)$ be the polynomial map such that $f^N\circ\sigma=\sigma\circ D$, with $\sigma$ the Poincaré map of $\sigma_0$ given by Proposition \ref{prop:thedevolepingmap}. We assume that the following formula holds  
\begin{equation}\label{eq:Dupont-Taflinformulabeforepullback}
    \mu=T\wedge dd^c|W_a|^2\ \mathrm{on}\ (\mathbb{P}^2,a).
\end{equation}
We also assume that there exists $\Omega$ a open set charged by $\mu$ such that for any open set $V\subset\Omega$, if $\mu(V)>0$ then $\mu^s(V)>0$. In Lemma \ref{rmq:remarqueDinh2}, we will in fact prove that it implies that, for any open set $V\subset\mathbb{P}^2$, if $\mu(V)>0$ then $\mu^s(V)>0$. We also assume that the exceptional set $\mathcal{E}$ of $f$ does not intersect $\mathrm{Supp}(\mu)$.

\subsection{Pull back by the Poincaré map}\label{sec:pullbackbythedevelopingmap}

Our purpose in this section is to see how (\ref{eq:Dupont-Taflinformulabeforepullback})  is lifted by  $\sigma$ on $\cmplex^2$. Let us specify how the positive closed $(1,1)-$current $\sigma^*T$ and the Borel measure $\sigma^*\mu$ are defined on $\cmplex^2$. We refer to {\cite{bedtay82, forsib94, sib99, dinsib}} for general accounts on currents.

For every $x \in \cmplex^2$, let us write $T=dd^c(u)$ on an open neighborhood of $\sigma(x)$, $u$ being a bounded continuous \emph{psh} function. Then $\sigma^*T$ is defined near $x$ by $dd^c (u \circ \sigma)$. Since $u$ is bounded, $\sigma^* \mu := \sigma^* T \wedge \sigma^* T$ is well defined near $x$ in the sense of Bedford-Taylor. These currents are representable by integration (they have order 0), hence the  trace of $\sigma^*T$ and $\sigma^*\mu$ are Borel measures on $\cmplex^2$. These two measures are finite on bounded open subsets by Chern-Levine-Nirenberg inequality, thus they are Borel regular measures (Radon measures) on $\cmplex^2$. 

\begin{prop}\label{prop:pullbackofduponttaflinformulabysigma} Let $f$ be a degree $d\geq2$ map on $\mathbb{P}^2$ of equilibrium measure $\mu=T\wedge T$. Assume the first hypothesis of $\mathrm{Theorem\ \ref{thm:Main}}$ so that the formula $(\ref{eq:Dupont-Taflinformulabeforepullback})$ holds for a repelling $N-$periodic point $a\in \mathrm{Supp}(\mu)$. With the preceding notations,
\begin{enumerate}
    \item $|\chi_2|=\sqrt{d^N}$.
    \item The pull back $\sigma^*\mu$ and the pull back $\sigma^*T$ are related on $\cmplex^2$ by
        $$\sigma^*\mu = \sigma^*T\wedge dd^c|w|^2\ \mathrm{on}\ \cmplex^2,$$
    where $w$ is the second standard coordinate on $\cmplex^2$.
    \item One has $\mu\ll\sigma_T$ on $\mathbb{P}^2$.
\end{enumerate}
\end{prop}
\noindent\textbf{\underline{Proof :}} Without loss of generality we can assume that $N=1$. Using $f^*\mu=d^2\mu$, $f^*T=dT$ and $W_a\circ f=\chi_2W_a$, we get by pulling back  (\ref{eq:Dupont-Taflinformulabeforepullback}) by $f$ :
$$\mu=\frac{|\chi_2|^2}{d}\left(T\wedge dd^c|W_a|^2\right) \textrm{ on } (\mathbb{P}^2,a) .$$
Using again (\ref{eq:Dupont-Taflinformulabeforepullback}) one gets $\mu=\frac{|\chi_2|^2}{d}\mu$ on $(\mathbb{P}^2,a)$. Finally $|\chi_2|^2=d$ since $a \in \mathrm{Supp}(\mu)$.

Let us prove the second item. We want to show that $\sigma^*T\wedge dd^c|w|^2$ and $\sigma^*\mu$ are equal. Since they are Borel regular measures, it suffices to prove that they coincide on compact sets $K$ of $\cmplex^2$. Let us fix a compact set $K\subset\cmplex^2$. Let $k$ be large enough such that $K \subset \mathbb{D}^2_k$, in particular $\sigma = \sigma_k = f^{n_k} \circ \sigma_0 \circ D^{-n_k}$ on $K$. Since ${f^{n_k}}^{*}\mu = d^{2n_k}\mu$, we get
\begin{align*}
    (\sigma^*\mu)(K) & = \left[(D^{-n_k})^*\sigma_0^*({f^{n_k}}^*\mu)\right](K)
                     = d^{2n_k}\times \left[ {(D^{-n_k})}^*\sigma_0^*\mu \right](K) .
\end{align*}
We obtain by the formula (\ref{eq:Dupont-Taflinformulabeforepullback}) :
\begin{align*}
    (\sigma^*\mu)(K) & = d^{2n_k}\times \left[(D^{-n_k})^*\sigma_0^*T)\wedge (dd^c|W_a\circ\sigma_0\circ D^{-n_k}|^2)\right](K).
\end{align*}
Now using $$W_a\circ\sigma_0\circ D^{-n_k}=w\circ\sigma_0^{-1}\circ\sigma_0\circ D^{-n_k}=w\circ D^{-n_k}=\chi_2^{-n_k}\times w , $$ we deduce that :
\begin{align*}
    (\sigma^*\mu)(K) & = d^{2n_k}|\chi_2|^{-2n_k}\times\left[(D^{-n_k})^*\sigma_0^*T)\wedge dd^c|w|^2\right](K).
\end{align*}
But  $|\chi_2|=\sqrt{d}$ by the first item, hence we obtain :
\begin{align*}
    (\sigma^*\mu)(K) & = [(D^{-n_k})^*\sigma_0^*({d^{n_k}T)})\wedge dd^c|w|^2](K) .
\end{align*}
Using the invariant relation $d^{n_k}T={f^{n_k}}^{*}T$ one deduces :
\begin{align*}
    (\sigma^*\mu)(K) & = \left[(D^{-n_k})^*\sigma_0^*{f^{n_k}}^*T)\wedge dd^c|w|^2\right](K).
\end{align*}
Finally, by recalling that $\sigma|_K=f^{n_k}\circ\sigma_0\circ D^{-n_k}|_K$, we get 
\begin{align*}
    (\sigma^*\mu)(K) & = \left[(\sigma^*T)\wedge dd^c|w|^2\right](K).
\end{align*}
This proves that the measures $\sigma^*\mu$ and $\sigma^*T\wedge dd^c|w|^2$ are equal on $\cmplex^2$.

Let us explain now why $\mu\ll\sigma_T$ on $\mathbb{P}^2$. Let $p$ be a point in $\mathrm{Supp}(\mu)\cap \sigma(\cmplex^2)$ which is not a critical point of $\sigma$. Then there exists a inverse branch $\sigma_p^{-1}$ of $\sigma$ defined on a neighborhood of $p$, and the formula $\sigma^*\mu=\sigma^*T\wedge dd^c|w|^2$ implies that we have $\mu=T\wedge dd^c|w\circ \sigma_p^{-1}|^2\ll\sigma_T$ near $p$. Since it is true for any point $p$ in $\mathrm{Supp}(\mu)\cap\sigma(\cmplex^2)\backslash \sigma(\mathrm{Crit}\ \sigma)$, we deduce that $\mu\ll \sigma_T$ on $\sigma(\cmplex^2)\backslash\sigma(\mathrm{Crit}\ \sigma)$. Recall that $\sigma(\cmplex^2)\supset\mathbb{P}^2\backslash\mathcal{E}$ (see Proposition \ref{prop:imageofdevelopmentmapssigmacontainsP^2minusE(f)}) and that $\mu$ does not charge the set $\mathcal{E}\cup \sigma(\mathrm{Crit}\ \sigma)$, see Proposition \ref{prop:thedevolepingmap}. We thus obtain $\mu\ll\sigma_T$ on $\mathbb{P}^2$.\qed 

\subsection{Proof of Theorem {\ref{thm:Main}}}\label{sec:finpreuvetheoremeB}

In the following remaining sections, we are going to prove the following theorem:

\begin{thm}\label{thm:foliatedJ'} Assume the hypothesis of $\mathrm{Theorem}\  \ref{thm:Main}$ with the same notations. Let $J$ be a compact subset of $\mathrm{Supp}(T)\backslash\mathcal{E}$ containing $\mathrm{Supp}(\mu)$ such that $J=\mathrm{Supp}(\sigma_T|_J)$. Then there exists a open neighborhood $\mathcal{V}$ of $J$, and a foliation $\mathcal{F}$ defined on $\mathcal{V}$, such that $\sigma^*\mathcal{F}=\mathcal{F}_w$ on $\sigma^{-1}(\mathcal{V})$, where $\mathcal{F}_w$ is the horizontal foliation of $\cmplex^2$.
\end{thm}

We then assert that this theorem implies Theorem \ref{thm:Main}, let us explain why. We start by explaining that $\mathrm{Supp}(\mu)$ is included in $\mathrm{Supp}(\mu^s)\backslash\mathcal{E}$. Since $\mathrm{Supp}(\mu)\cap\mathcal{E}=\emptyset$, it remains to prove:
\begin{lemme}\label{rmq:remarqueDinh2}
    Item \textit{2.} of $\mathrm{Theorem\ \ref{thm:Main}}$ is equivalent to the following:
    \begin{enumerate}
    \item[2'.] Any open set charged by $\mu$ is also charged by $\mu^s$.
    \end{enumerate}
    In particular $\mathrm{Supp}(\mu)\subset\mathrm{Supp}(\mu^s)$ and thus $\mathrm{Supp}(\mu^s)=\mathrm{Supp}(T)$.
\end{lemme}
\noindent\textbf{\underline{Proof:}} Indeed, this assumption implies \textit{2.} Let us give briefly the arguments for the reverse implication. Let $\Omega$ be the open neighborhood provided by item \textit{2.} in the Theorem \ref{thm:Main}. By assumption, $\mu^s$ charges every open subset of $\Omega$ which is charged by $\mu$. Let $V$ be a neighborhood of a point in $\mathrm{Supp}(\mu)$. Since $\mu$ is mixing there exists $n_0$ such that $\mu(f^{-n_0}\Omega\ \cap V)>0$. Let $x_0\in\left(\mathrm{Supp}(\mu)\cap f^{-n_0}\Omega\ \cap V\right)\backslash \mathrm{Crit}(f^{n_0})$ and let $B \subset V$ be a small ball centered at $x_0$ such that $f^{n_0}(B) \subset\Omega$ and  $f^{n_0}$ is injective on $B$. Observe that $\mu^s(f^{n_0}(B))>0$ by hypothesis. Since $(f^{n_0})^*T = d^{n_0}T$ and since $\omega_{\mathbb{P}^2}$ is quasi-invariant by the biholomorphism $f^{n_0} : B \to f^{n_0}(B)$, the measures $\sigma_T\circ f^{n_0}$ and $\sigma_T$ are equivalent on $B$. Moreover the measures $\mu \circ f^{n_0}$ and $\mu$ are also equivalent on $B$ since $(f^{n_0})^*\mu = d^{2n_0} \mu$. Hence $\mu^s(B) > 0$ as desired.\qed\\

Consider now a covering of $\mathbb{P}^2\backslash\mathcal{E}$ by compact sets
$\mathbb{P}^2\backslash\mathcal{E}=\bigcup_{j=1}^{+\infty}K_j,$
where the $K_j's$ satisfy $K_j\subset \overset{\circ}{K}_{j+1}$. Then intersecting this covering with $\mathrm{Supp}(T)$ we obtain a covering
$\mathrm{Supp}(T)\backslash\mathcal{E}=\bigcup_{j=1}^{+\infty}L_j,$
where $L_j:=\mathrm{Supp}(T)\cap K_j$. Up to re-index the sequence $(K_j)_{j\geq1}$, we can assume that for any $j\geq 1$:
\begin{center}
    $\mathrm{Supp}(\mu)\subset \overset{\circ}{K}_j\cap\mathrm{Supp}(T)\subset L_j\subset\mathrm{Supp}(T)\backslash\mathcal{E}$.
\end{center}
Observe that these inclusions are possible since $\mathrm{Supp}(\mu)\subset\mathrm{Supp}(T)$ (by Proposition \ref{prop:pullbackofduponttaflinformulabysigma} or by Lemma \ref{rmq:remarqueDinh2}),
and since $\mathrm{Supp}(\mu)\cap\mathcal{E}=\emptyset$.
\begin{lemme}\label{lemma:lemmesupplementairesurlessupports}
    Let $J_j:=\mathrm{Supp}(\sigma_T|_{L_j})$ for any $j\geq 1$. Then $(J_j)_{j\geq1}$ is a increasing sequence of compact sets such that for any $j\geq1$:
    \begin{enumerate}
        \item $J_j=\mathrm{Supp}(\sigma_T|_{J_j})$ and $\mathrm{Supp}(\mu)\subset J_j\subset \mathrm{Supp}(T)\backslash\mathcal{E}$.
        \item $\mathrm{Supp}(T)\backslash\mathcal{E}=\bigcup_{j=1}^{+\infty}J_j$.
    \end{enumerate}
\end{lemme}
\noindent\underline{\textbf{Proof:}} Recall that $\mathrm{Supp}(T)=\mathrm{Supp}(\sigma_T)$ and that $\mu\ll\sigma_T$ by Proposition \ref{prop:pullbackofduponttaflinformulabysigma}.\\ 
\textit{1.} First $\mathrm{Supp}(\sigma_T|_{J_j})\subset J_j$, second if $p\in J_j$ then $\sigma_T(J_j\cap B(p,\varepsilon))=\sigma_T(L_j\cap B(p,\varepsilon))$ by definition of $J_j$, but $p\in J_j=\mathrm{Supp(\sigma_T|_{L_j})}$ thus $\sigma_T(L_j\cap B(p,\varepsilon))>0$, for any $\varepsilon>0$. We deduce $\mathrm{Supp}(\sigma_T|_{J_j})=J_j$. To conclude, observe that because $\mu\ll\sigma_T$ and $\mathrm{Supp}(\mu)\subset \overset{\circ}{K}_j$, we must have $\mathrm{Supp}(\mu)\subset\mathrm{Supp}(\sigma_T|_{L_j})=J_j$. 

\noindent \textit{2.} Let $p\in\mathrm{Supp}(T)\backslash\mathcal{E}$, then there exists $j\geq1$ such that $p\in L_j$, let us assume $j=1$ for simplicity. Then in particular $p\in \overset{\circ}{K}_2$ and there exists $\varepsilon_0$ such that $B(p,\varepsilon_0)\subset K_2$. Then for every $0<\varepsilon\leq \varepsilon_0$, we have $\sigma_T(L_2\cap B(p,\varepsilon))=\sigma_T(\mathrm{Supp}(T)\cap K_2\cap B(p,\varepsilon)) = \sigma_T(K_2\cap B(p,\varepsilon))$ by definition of $L_2$. We deduce $\sigma_T(L_2\cap B(p,\varepsilon))=\sigma_T(B(p,\varepsilon))>0$ since $B(p,\varepsilon)\subset K_2$ and since $p\in L_1\subset \mathrm{Supp}(T)$. It implies $p\in\mathrm{Supp}(\sigma_T|_{L_2})=J_2$ and the conclusion follows.\qed\\

Theorem \ref{thm:foliatedJ'} and Lemma \ref{lemma:lemmesupplementairesurlessupports} thus imply the existence of a sequence of open sets $(\mathcal{V}_j)_{j\geq1}$ and of a sequence of foliations $(\mathcal{F}_j)_{j\geq1}$ such that for any $j\geq1$:
\begin{center}
    $J_j\subset \mathcal{V}_j$ and $\mathcal{F}_j$ is defined on $\mathcal{V}_j$, and $\sigma^*\mathcal{F}_j=\mathcal{F}_w$ on $\sigma^{-1}(\mathcal{V}_j)$.
\end{center}
\begin{lemme}\label{lemma:F=F_1surV}
    For any $j\geq2$, $\mathcal{F}_j=\mathcal{F}_{j-1}$ on $\mathcal{V}_j\cap \mathcal{V}_{j-1}$.
\end{lemme}
\noindent\textbf{\underline{Proof:}} Let us explain the argument for $j=2$. We have
$\sigma^*\mathcal{F}_2=\mathcal{F}_{w}=\sigma^*\mathcal{F}_1$ on $\sigma^{-1}(\mathcal{V}_2\cap \mathcal{V}_1)$. Let $\mathcal{C}$ be a connected component of $\mathcal{V}_2\cap \mathcal{V}_1$, and let $U\neq\emptyset$ be an open set such that $\sigma(U)\subset \mathcal{C}$ and such that $\sigma$ is injective on $U$. Then observe that $\mathcal{F}_2=\mathcal{F}_1$ on $\sigma(U)$ by applying $(\sigma|_U)_*$ to the equality $\sigma^*\mathcal{F}_2|_U=\sigma^*\mathcal{F}_1|_U$. It implies $\mathcal{F}_2=\mathcal{F}_1$ on $\mathcal{C}$ by analytic continuation (cf. Lemma \ref{lemma:F_1=F_2onD^2(1)thenF_1=F_2onD^2(2)}), for every $\mathcal{C}$. The result follows.\qed\\ 

This lemma implies that the foliations $\mathcal{F}_j,\ j\geq 1$, patch together and create a foliation 
$$\mathcal{F}:=\bigcup_{j=1}^{+\infty}\mathcal{F}_j$$
on the open set $\mathcal{V}:=\left(\bigcup_{j=1}^{+\infty}\mathcal{V}_j\right) \supset \mathrm{Supp}(T)\backslash\mathcal{E}$. Observe that we have $\sigma^*\mathcal{F}=\bigcup_{j=1}^{+\infty}\sigma^*\mathcal{F}_j=\mathcal{F}_{w}\ \mathrm{on}\ \sigma^{-1}(\mathcal{V}).$
So we can conclude that $\mathcal{F}$ is invariant on $\mathcal{V}$ by $f^N$ (recall that $N$ is the period of the repelling point $a$):

\begin{prop}\label{prop:invarianceofFbyfN}
    $(f^N)^*\mathcal{F}=\mathcal{F}$ on $f^{-N}(\mathcal{V})\cap\mathcal{V}$.
\end{prop}
\noindent\textbf{\underline{Proof :}} We can assume $N=1$ without loss of generality. Let us fix arbitrary points $p,q\in \mathcal{V}$ such that $f(p)=q$. Let us prove that there exists $W\subset\mathcal{V}$ an open neighborhood of $p$, and $V\subset\mathcal{V}$ an open neighborhood of $q$, satisfying $f(W)\subset V$, and such that $(f^{*}\mathcal{F}|_V)|_W=\mathcal{F}|_W$. It is enough to ensure the result.

Let $x\in\sigma^{-1}(p)$ and let $y:=D(x)$ which satisfies $\sigma(y)=q$. We can consider a connected open set $U_x$ (resp. $U_y$) in $\sigma^{-1}(\mathcal{V})$ containing $x$ (resp. $y$) such that $W:=\sigma(U_x)$ (resp. $V:=\sigma(U_y)$) is contained in $\mathcal{V}$. We can assume $D(U_x)\subset U_y$ and $f(W)\subset V$. We denote $\sigma_p:=\sigma|_{U_x}$ and $\sigma_q:=\sigma|_{U_y}$.

Let now $\mathcal{G}:=(f^*\mathcal{F}|_V)|_W$, observe that $\sigma_p^*\mathcal{G}=(D^*\sigma_q^*\mathcal{F})|_{U_x}$ by using that $\sigma\circ D=\sigma_q\circ D$ on $U_x\subset D^{-1}(U_y)$. Since $\sigma^*\mathcal{F}=\mathcal{F}_w$ on $\sigma^{-1}(\mathcal{V})$, we have $\sigma_q^*\mathcal{F}=\mathcal{F}_w|_{U_y}$ and thus $\sigma_p^*\mathcal{G}=(D^*\mathcal{F}_w)|_{U_x}=\mathcal{F}_w|_{U_x}$. The equality $\sigma^*\mathcal{F}=\mathcal{F}_w$ on $\sigma^{-1}(\mathcal{V})$, also implies $\mathcal{F}_w|_{U_x}=\sigma_p^*(\mathcal{F}|_W)$, and we deduce that $\sigma_p^*\mathcal{G}=\sigma_p^*(\mathcal{F}|_W)$. It implies that $\mathcal{G}=\mathcal{F}|_{W}$ since $U_x$ is connected (pull-back by a local inverse of $\sigma$ on $U_x$ and conclude by analytic continuation, cf. Lemma \ref{lemma:F_1=F_2onD^2(1)thenF_1=F_2onD^2(2)}).\qed

\subsection{Proof of Theorem \ref{thm:foliatedJ'}}\label{sec:proofoftheorem62}

Let us fixed from now on a compact set $J\subset \mathrm{Supp}(T)\backslash\mathcal{E}$ such that $J\supset \mathrm{Supp}(\mu)$, and such that $\mathrm{Supp}(\sigma_T|_J)=J$. Our goal is now to prove that there exists an open neighborhood $\mathcal{V}$ of $J$ and a foliation $\mathcal{F}$ defined on $\mathcal{V}$ which satisfies $\sigma^{*}\mathcal{F}=\mathcal{F}_w$ on $\mathcal{V}$. We proceed in several steps.

\subsubsection{Step 1: Inverse branches near regular values of \texorpdfstring{$J$}{TEXT}}\label{sec:inversebranchesaboveregularvaluesofsigmainsidethejuliaset}

Recall that $J\supset \mathrm{Supp}(\mu)$ and that $J\cap\mathcal{E}=\emptyset$. Thus we have that $J \subset \sigma(\cmplex^2)$ by Proposition \ref{prop:imageofdevelopmentmapssigmacontainsP^2minusE(f)}. 

By Proposition \ref{prop:thedevolepingmap}, the map $\sigma$ is equal to
$\sigma_k=f^{n_k}\circ\sigma_0\circ D^{-n_k}$ on $\mathbb{D}^2_k$. Let now $p\in J \subset\sigma(\cmplex^2)$ and let $k \geq 1$ be such that $p \in \sigma(\mathbb{D}_k^2)$. Since $\sigma$ is open, $\sigma(\mathbb{D}^2_k)$ is an open neighborhood of $p$. By compactness we can cover $J$  by a finite number of such neighborhoods
$$J\subset \mathcal{U} := \sigma(\mathbb{D}^2_{k_1})\cup\cdots\cup\sigma(\mathbb{D}^2_{k_N}).$$
Let $R:=\underset{1\leq i\leq N}{\max}k_i$ and observe that $\underset{j=1}{\overset{N}{\cup}}\mathbb{D}^2_{{k_j}}= \mathbb{D}^2_R=:D_R$, thus we have a map:
\begin{equation*}
    \label{eq:sigma(D_R)supsetU}
    \sigma_R :=(f^{n_R} \circ\sigma_0\circ D^{-n_R}) |_{D_R} : D_R\longrightarrow \sigma_R (D_R)= {\mathcal{U}}\supset J.
\end{equation*}
According to Proposition \ref{prop:thedevolepingmap}, the critical values $E_R$ of $\sigma_R$ satisfies:
$$E_R=\sigma_R(\mathrm{Crit}\ \sigma_R) \subset f^{n_R}(\mathrm{Crit}\ f^{n_R}).$$
In particular $\mathrm{Reg}_{\sigma_R}(J):=J\backslash E_R$ has full $\sigma_T|_J-$measure and full $\mu-$measure.

Any point $p\in \mathcal{U}$ satisfies $1\leq \#(\sigma_R^{-1}(p))\leq d^{2 n_R}$. So for every  $p\in\mathrm{Reg}_{\sigma_R}(J)$ there exists $n_p\in\{1,\cdots,d^{2 n_R}\}$ such that $\sigma^{-1}_R(p)=\{x_{1,p},\cdots,x_{n_p,p}\}$. For every $j \in \{ 1, \ldots , n_p \}$ there exists  an open connected neighborhood $U_{j,p} \subset D_R$ of $x_{j,p}$ such that $$\sigma_{j,p}:=\sigma_R |_{U_{j,p}}:U_{j,p}\longrightarrow\sigma(U_{j,p})=:V_{j,p} \subset \mathcal{U}$$ is a biholomorphism. We denote $\sigma_{j,p}^{-1}:V_{j,p}\to U_{j,p}$ the inverse map. We can assume that every $V_{j,p}$ does not intersect $E_R$. We denote 
\begin{equation*}
    V_p:=\bigcap_{j=1}^{n_p}V_{j,p}.
\end{equation*}
We can assume that $V_p$ is connected. 

Observe that, even up to a reduction of $V_p$, the fiber of a point $p'\in V_p$ by $\sigma_R$ may not be given by the inverse branches $\sigma^{-1}_{1,p},\cdots,\sigma_{n_p,p}^{-1}$. Indeed, if $p$ admits a preimage by $\sigma$ on the boundary of $D_R$, the inclusion $\left\{\sigma^{-1}_{1,p}(p'),\cdots,\sigma_{n_p,p}^{-1}(p')\right\}  \subset \sigma_R^{-1}(p')$ could be strict. To avoid this difficulty, we shall introduce special sets of inverse branches :

\begin{defn}\label{defn:InverseBranchesIB_R(p)} For every $p\in\mathrm{Reg}_{\sigma_R}(J)$, we define 
\begin{equation}\label{eq:localexpressionofthefibersofsigma_Rbythesectionsabovetheregularpointp}
    \mathrm{IB}_R(p):=\bigcup_{j=1}^{n_p}\sigma_{j,p}^{-1}(V_p) \subset D_R.
\end{equation}
Here the abreviation "$\mathrm{IB}$" stands for \textit{inverse branches}.
\end{defn}

\subsubsection{Step 2: Construction of the foliation near regular values}\label{sec:constructionofthefoliationnearregularpoints}

In this section we introduce for every $p\in\mathrm{Reg}_{\sigma_R}(J)$ a foliation $\mathcal{F}_p$ on $V_p$. We shall patch together the foliations $(\mathcal{F}_p)_{p\in\mathrm{Reg}_{\sigma_R}(J)}$ (and also the foliations on neighborhoods of singular values $q\in J\cap E_R$ constructed in Section \ref{sec:donctructionofthefoliationnearsingularvalues}). To do so, we shall use the following distribution of tangent complex lines, defined for each for $p\in\mathrm{Reg}_{\sigma_R}(J)$:
\begin{equation}\label{eq:Directionoffoliationmap}
\mathcal{D}_p : 
\left\{
    \begin{array}{ll}
      V_p & \longrightarrow\ \mathbb{P}(T\mathbb{P}^2)\\
      \textcolor{white}{0}      & \textcolor{white}{0}\\
      p'                         & \longmapsto\ \left[d_{x'}\sigma_R \cdot(1,0)\right]\ \mathrm{for\ any}\ x'\in \mathrm{IB}_R(p)\ \mathrm{such\ that}\ \sigma_R(x')=p'\\
    \end{array}
\right.
\end{equation}
It will describe the directions followed by the different local foliations we will construct. The notation $[\vec v]$ stands for the complex line $\cmplex\cdot\vec v$ of $T_p\mathbb{P}^2$. We verify in Proposition \ref{prop:generalcaseoffiberconnections} that the maps $\mathcal{D}_p$ are well defined. The proof of this proposition uses the Patching Theorem \ref{lemme:patchinglemma}, so wee need that $\mu^s$ charges all the neighborhoods of points of $\mathrm{Supp}(\mu)$. This property is ensured by Item \textit{2.} of Theorem \ref{thm:Main} and Lemma {\ref{rmq:remarqueDinh2}}.

\begin{prop}\label{prop:generalcaseoffiberconnections} Assume that $f$ and the $N-$periodic repelling point $a$ satisfy the hypothesis of $\mathrm{Theorem\ \ref{thm:Main}}$ with the same notations. 
\begin{enumerate}
    \item Let $x,y\in \cmplex^{2}\backslash(\mathrm{Crit}\ \sigma)$ and let $p:=\sigma(x)$ and $q:=\sigma(y)$. Let $U_{x,p}$ be an open connected neighborhood of $x$ such that the map $\sigma_{x,p}:=\sigma|_{U_{x,p}}:U_{x,p}{\to} \sigma(U_{x,p})=:V_{x,p}$ is a biholomorphism. Similarly let $\sigma_{y,q}:=\sigma|_{U_{y,q}}:U_{y,q}\to\sigma(U_{y,q})=:V_{y,q}$ be a biholomorphism. Assume that $\Omega:=V_{x,p}\cap V_{y,q}$ is connected and that $\sigma_T(\Omega)>0$. Then $\phi:=\left(\sigma^{-1}_{y,q}\circ\sigma_{x,p}\right)$ has the form $\phi(z,w)=(A(z,w),B(w))$ on $\sigma_{x,p}^{-1}(\Omega)$.
    \item Let $p\in \mathrm{Reg}_{\sigma_R}(J)$ and $p'\in V_p$. Then for every $x',y'\in\sigma_R^{-1}(p')\cap\mathrm{IB}_R(p)$, the complex lines $[d_{x'}\sigma\cdot(1,0)]$ and $[d_{y'}\sigma\cdot(1,0)]$ are equal. In particular the map $\mathcal{D}_p: V_p \to\ \mathbb{P}(T\mathbb{P}^2)$ is well defined.
    \item Let $p,q\in\mathrm{Reg}_{\sigma_R}(J)$. If $\mathcal{C}$ is a connected component of $V_p\cap V_q$ of measure $\sigma_T(\mathcal{C})>0$, then $\mathcal{D}_p=\mathcal{D}_q$ on $\mathcal{C}$.
\end{enumerate}
\end{prop}
\noindent\textbf{\underline{Proof :}} \ \\
\noindent\textit{1.} Let us write $\phi=(A(z,w),B(z,w))$ and let us prove that $B(z,w)$ does not depend on $z$. To do so we use that $\sigma^*\mu = \sigma ^*T\wedge dd^c|w|^2\ \mathrm{on}\ \cmplex^2$, according to Proposition \ref{prop:pullbackofduponttaflinformulabysigma}. Restricting this equality on $U_{x,p}$ (resp. on $U_{y,q}$) and pushing forward it by $\sigma_{x,p}$ (resp. by $\sigma_{y,q}$) we get on $\Omega=V_{x,p}\cap V_{y,q}$ :
$$\mu=T\wedge dd^c|W_p|^2\ \mathrm{and}\ \mu=T\wedge dd^c|W_q|^2,$$
with the notations $(Z_p,W_p):=\sigma_{x,p}^{-1}\ \mathrm{and}\ (Z_q,W_q):=\sigma_{y,q}^{-1}$.

We recall that $\sigma_T=\mu^a + \mu^s$ is defined by \eqref{eq:decompositionofRNofsigma_Twithrespecttomu}. 
By assumption $\sigma_T(\Omega)>0$ so by Lemma \ref{rmq:remarqueDinh2} it ensures that 
$\mu^s(\Omega)>0$. By hypothesis, $\Omega$ is also connected, thus we can apply the Patching Theorem \ref{lemme:patchinglemma}: there exists $\beta\in\mathcal{O}^*(\Omega)$ such that $dW_q=\beta dW_p$ which implies $dB=\sigma_{x,p}^*dW_q=(\beta\circ\sigma_{x,p})dw$ on $\sigma_{x,p}^{-1}(\Omega)$. In particular $\partial_zB\equiv0$ on $\sigma_{x,p}^{-1}(\Omega)$.

\noindent\textit{2.} Since $x',y'\in\mathrm{IB}_R(p)$ and  $\sigma(x')=\sigma(y')=p'$, there exist $i,j$ such that $x'=\sigma_{i,p}^{-1}(p')$ and $y'=\sigma_{j,p}^{-1}(p')$. Let us assume $i=1$ and $j=2$ for simplicity. Observe that  $\sigma_{1,p}^{-1}$ and $\sigma_{2,p}^{-1}$ are defined on the connected set $V_p$ by construction, and $\sigma_T(V_p)>0$ since $V_p$ is an open neighborhood of $p \in J\subset\mathrm{Supp}(T)$. So we can apply the point \textit{1.} with $p=q$, to conclude that the map $\phi := \left(\sigma_{2,p}^{-1}\circ \sigma_{1,p}\right)$ satisfies $\phi(z,w)=(A(z,w),B(w))$ on $\sigma_{1,p}^{-1}(V_p)$.
In particular we have $d_{x'}\phi\cdot(1,0)=\partial_z A(x')\cdot(1,0)$ and so we deduce:
$$d_{x'}\sigma_{1,p}\cdot(1,0)=d_{y'}\sigma_{2,p}\cdot d_{x'}\phi\cdot(1,0)=\partial_z A(x')\times d_{y'}\sigma_{2,p}\cdot(1,0).$$
The coefficient $\partial_zA(x')$ is not equal to $0$ because $\phi$ is a biholomorphism. Finally
$$\left[d_{x'}\sigma\cdot(1,0)\right] = \left[d_{y'}\sigma\cdot(1,0)\right]\in \mathbb{P}\left(T_{p'}\mathbb{P}^2\right),$$
which proves that $\mathcal{D}_p$ is well defined on $V_p$.

\noindent\textit{3.} Let $\mathcal{C}$ be a connected component of $V_p\cap V_q$, and let us assume that $\sigma_T(\mathcal{C})>0$. Let us fix $p'$ an element of $\mathcal{C}$, and let $x\in\sigma_R^{-1}(p')\cap\mathrm{IB}_R(p)$ and $y\in\sigma_R^{-1}(p')\cap\mathrm{IB}_R(q)$. By definition of $\mathrm{IB}_R(p)$ there exists $i\in\{1,\cdots,n_p\}$ such that $x=\sigma_{i,p}^{-1}(p')$. Similarly there exists $j\in\{1,\cdots,n_q\}$ such that $y=\sigma_{j,q}^{-1}(p')$. Let now $U_x:=\sigma_{i,p}^{-1}(\mathcal{C})\subset U_{i,p}$ be a neighborhood of $x$, $\mathcal{C}\subset V_p$ is indeed contained in the domain of definition of $\sigma_{i,p}^{-1}$. Similarly let us define $U_y:=\sigma_{j,q}^{-1}(\mathcal{C})\subset U_{j,q}$ which contains $y$. We have $\sigma(U_x) = \sigma(U_y) =\mathcal{C}$, and by hypothesis $\mathcal{C}$ is a connected open set such that $\sigma_T(\mathcal{C})>0$. The map $\phi:=(\sigma_{j,q}^{-1}\circ \sigma_{i,p}):U_x\to U_y$ is then well defined, and according to the first item (applied with $\Omega=\mathcal{C}$), $\phi$ has the form $\phi(z,w)=(A(z,w),B(w))$ on $U_x$. As in the proof of the preceding item, we can compute the matrix $d_x\phi$ and check that $d_x\sigma\cdot(1,0)=\partial_zA(x)\times d_y\sigma\cdot(1,0)$. We also have $\partial_zA(x)\neq0$ and thus using again the second item we have:
$$\mathcal{D}_p(p') = [d_{x}\sigma\cdot(1,0)] = [d_{y}\sigma\cdot(1,0)] = \mathcal{D}_q(p').$$
The proof is then complete.\qed 

\begin{lemme}\label{lemma:omega_j,p=d(wcircsigma_j,p^-1)} Let $p\in\mathrm{Reg}_{\sigma_R}(J)$. Then the holomorphic 1-forms  $$\omega_{j,p}:=d(w\circ\sigma^{-1}_{j,p}),\ j\in\{1,\cdots,n_p\},$$
are equal modulo a multiplicative function in $\mathcal{O}^*(V_p)$. 
\end{lemme}

\noindent\textbf{\underline{Proof :}} By the point \textit{1.} of Proposition \ref{prop:generalcaseoffiberconnections}, for each $i,j$ we have on $V_p$ that $\sigma^{-1}_{j,p}\circ\sigma_{i,p}=(A_{ij}(z,w),B_{ij}(w))$
and thus $\sigma_{i,p}^*\omega_{j,p} = d(w\circ\sigma^{-1}_{j,p}\circ\sigma_{i,p}) = B_{ij}'(w)dw$. Then $\omega_{j,p} = (B_{ij}'\circ (w \circ \sigma^{-1}_{i,p}))\times \omega_{i,p}$. Because the function $w\circ\sigma_{i,p}^{-1}$ is a submersion, we have $(B_{ij}'\circ(w \circ \sigma^{-1}_{i,p}))\in\mathcal{O}^*(V_p)$. \qed \\

Thanks to Lemma \ref{lemma:omega_j,p=d(wcircsigma_j,p^-1)} we can put the following definition. We refer to Section \ref{sec:holomorphicfoliations} for the definition of foliations using $1$-forms. Recall that $\mathcal{F}_w$ is the horizontal foliation on $\mathbb{C}^2$ given by the $1$-form $dw$. We also recall that $\mathrm{IB}_R(p)$ is defined by \eqref{eq:localexpressionofthefibersofsigma_Rbythesectionsabovetheregularpointp}.

\begin{defn}\label{defn:definitionofF_pbyomega_pandsigma^*F=F_wonIB_R(p)}
For every $p\in\mathrm{Reg}_{\sigma_R}(J)$ we denote $\mathcal{F}_p$ the non singular foliation on $V_p$ defined by one of the holomorphic 1-forms $\omega_{j,p}$. We note that ${\sigma|_{\mathrm{IB}_R(p)}} ^*\mathcal{F}_p=\mathcal{F}_w$.
\end{defn}

Finally one can reformulate the results of the present section as follows.
 
\begin{prop}\label{prop:T_p'P^2=D(p')} \
\begin{enumerate}
    \item If $p\in\mathrm{Reg}_{\sigma_R}(J)$, then  $T_{p'}\mathcal{F}_p=\mathcal{D}_p(p')$ for every $p'\in V_p$.
    \item $\forall p,q\in\mathrm{Reg}_{\sigma_R}(J)$,  $\mathcal{F}_p$ and $\mathcal{F}_q$ coincide on every connected component $\mathcal{C}\subset V_p\cap V_q$ such that $\sigma_T(\mathcal{C})>0$
    (use \textit{1.}, $\mathcal{D}_p|_{\mathcal{C}}=\mathcal{D}_q|_{\mathcal{C}}$ by $\mathrm{Prop.}\ \ref{prop:generalcaseoffiberconnections}$ and $\mathrm{Lem.\ \ref{lemma:F_1=F_2onD^2(1)thenF_1=F_2onD^2(2)}}$). 
\end{enumerate}
\end{prop}

\subsubsection{Step 3: Construction of the foliation near singular values}\label{sec:donctructionofthefoliationnearsingularvalues}

For every singular value $q$ of $\sigma_R$ which belongs to $J$, we construct an open neighborhood $W_q\ni q$ and a foliation $\mathcal{F}_q$ (possibly with singularities) on $W_q$ which is tangent to directions $(\mathcal{D}_p)_{p\in\mathrm{Reg}_{\sigma_R}(J)}$ on neighborhoods of points of $\mathrm{Reg}_{\sigma_R}(J)\cap W_q$, see Proposition \ref{lemme:v_q(t)inD(t)}. Thanks to this property $\mathcal{F}_q$ will coincide with the regular foliations $(\mathcal{F}_{p})_{p\in\mathrm{Reg}_{\sigma_R}(J)}$ on connected components of $\bigcup_{p\in\mathrm{Reg}_{\sigma_R}(J)}V_p\cap W_q$ which have positive measure by $\sigma_T|_J$, see Proposition \ref{prop:patchingofFregwithF_qandF_qwithF_q'}. We finish the construction of the desired foliation in Section \ref{sec:finitecoveringofJ} using the compactness of $J$. 

To create $\mathcal{F}_q$ the idea is to construct, using $d\sigma\cdot(1,0)$, a vector field on a neighborhood of $q$. Let us fix a singular value $q\in J\cap E_R$. We recall that $\sigma(D_R)\supset J$ and according to Proposition \ref{prop:thedevolepingmap}, the fiber $\sigma^{-1}(q)$ is discrete in $\cmplex^2$. Thus we can consider a preimage $y\in D_R$ and an open set $U_y\subset \overline{U}_y\subset D_R$ containing $y$ such that $\overline{U}_y\cap \sigma^{-1}(q)=\{y\}$. Let $B_y$ be a centered ball at $y$ such that $\overline{B}_y\subset U_y$.

\begin{lemme}\label{lemma:constructionofW_yandW_qstsigma(W_y)=W_q} \ 
\begin{enumerate}
    \item There exists a connected neighborhood $U_{y,q}\subset B_y$ of $y$, such that $W_q:=\sigma(U_{y,q})$ is a connected neighborhood of $q$ satisfying: 
    \begin{equation}\label{eq:boitesdesecuritepouryetq}
        \sigma^{-1}(p)\cap U_y\subset B_y,\ \forall p\in W_q.
    \end{equation}
    \item Let $p\in W_q\backslash E_R$ and let $n_{p,q}$ be the cardinality of $\sigma^{-1}(p)\cap U_y$, it satisfies $1\leq n_{p,q}<+\infty$. There exist a ball $B_p\subset W_q\backslash E_R$ centered at $p$ and a family of inverse branches $\sigma_{1,p,q}^{-1},\cdots,\sigma_{n_{p,q},p,q}^{-1}$ such that for all $p'\in B_p$:
    \begin{equation}\label{eq:sigma^-1(p')capU_y=sigma_j^-1(p')j=1n_pq}
        \sigma^{-1}(p')\cap U_y = \left\{\sigma_{1,p,q}^{-1}(p'),\cdots,\sigma_{n_{p,q},p,q}^{-1}(p')\right\}.
    \end{equation}
\end{enumerate}
\end{lemme}
\noindent\textbf{\underline{Proof :}} \ \\
\textit{1.} Let us assume to the contrary that for any neighborhood of $q$ there exists a preimage by $\sigma$ of a point in this neighborhood which belongs to $U_{y}\backslash B_y$. Then looking at smaller and smaller neighborhoods of $q$ we obtain a sequence $(p_n)_n$ of points converging to $q$, and a sequence $(x_n)_n$ of $U_{y}\backslash B_y$ such that $\sigma(x_n)=p_n$. Taking $x\in \overline{U}_{y}\backslash B_y$ a cluster value of $(x_n)_n$, we have by continuity of $\sigma$ that, up to a sub-sequence, $\sigma(x)=\lim_{n} \sigma(x_{n})=\lim_{n} p_{n}=q$. Thus we have $x\in\sigma^{-1}(q)\cap \overline{U}_{y}\backslash B_y$, but this is a contradiction since $\sigma^{-1}(q)\cap \overline{U}_{y}=\{y\}\subset B_y$. 

So we have proved the existence of a neighborhood $W'_q$ of $q$ such that for all $p\in W'_q$, $\sigma^{-1}(p)\cap U_y\subset B_y$. Taking a small connected neighborhood $U_{y,q}$ of $y$ included in $B_y$, one has $\sigma(U_{y,q})\subset W_q'$. We complete the proof by setting $W_q:=\sigma(U_{y,q})$.\\

\noindent\textit{2.} Observe first that because $\sigma(B_y)\supset W_q$, the fiber $\sigma^{-1}(p)\cap U_y$ is not empty. As explained above the fibers of $\sigma$ are discrete in $\cmplex^2$, thus the fiber $\sigma^{-1}(p)\cap U_y$ is finite. In particular we have $1\leq n_{p,q}<+\infty$. Let $(\sigma_{j,p,q})_{j=1,\cdots,n_{p,q}}$ be a family of inverse branches of $\sigma$ such that $\sigma^{-1}(p)\cap U_y=\{\sigma_{j,p,q}^{-1}(p),\ 1\leq j\leq n_{p,q}\}$. Let us fix a ball $B_p\subset W_q\backslash E_R$ centered at $p$ such that all these inverse branches are defined on $B_p$. We denote $\mathrm{IB}(p,q):=\bigcup_{j=1}^{n_{p,q}}\sigma_{j,p,q}^{-1}(B_p)$. Since $\sigma^{-1}(p)\cap U_y\subset B_y$ by (\ref{eq:boitesdesecuritepouryetq}), up to a reduction of the radius of the ball $B_p$, we can assume that $\mathrm{IB}(p,q)\subset B_y$. Observe that $\mathrm{IB}(p,q)\subset\sigma^{-1}(B_p)\cap U_y$ by construction. Using \eqref{eq:boitesdesecuritepouryetq} and similar arguments involved in the previous item, we also have $\sigma^{-1}(B_p)\cap U_y\subset \mathrm{IB}(p,q)$, up to reduce $B_p$. Finally, observe that $\mathrm{IB}(p,q)=\sigma^{-1}(B_p)\cap U_y$ is exactly \eqref{eq:sigma^-1(p')capU_y=sigma_j^-1(p')j=1n_pq}.\qed 

\begin{defn}\label{defn:originaldefinitionofv_q} Let us fix $q\in J\cap E_R$. Let $W_q$ and $U_y$ be the open sets given by $\mathrm{Lemma\ \ref{lemma:constructionofW_yandW_qstsigma(W_y)=W_q}}$. We define a vector field $v_q:W_q\backslash E_R\longrightarrow T\mathbb{P}^2$ by:
$$\forall p\in W_q\backslash E_R,\ v_q(p):=\sum_{x\in\sigma^{-1}(p)\cap U_y}h(x)d_x\sigma\cdot(1,0)\in T_p\mathbb{P}^2.$$
The function $h\in\mathcal{O}(\cmplex^2)$ is chosen to have $v_q\not\equiv0$ on $W_q\backslash E_R$.
\end{defn}

\begin{lemme}\label{rmq:remarkonthevectorfieldvandthedefinitionofthevectorfieldsv_p} \
\begin{enumerate}
    \item The vector field $v_q:W_q\backslash E_R\longrightarrow T\mathbb{P}^2$ is holomorphic. 
    \item The function $h\in\mathcal{O}(\cmplex^2)$ chosen such that $v_q\not\equiv0$ exists. Moreover one can choose $h$ equal to a polynomial function on $\cmplex^2$.
\end{enumerate}
\end{lemme}
\noindent\textbf{\underline{Proof :}} \ \\
\noindent\textit{1.} Let $p\in W_q\backslash E_R$ fixed. Let $(\sigma_{j,p,q}^{-1})_{j\in\{1,\cdots,n_{p,q}\}}$ be the family of inverse branches defined on a ball $B_p$ given by the second item of Lemma \ref{lemma:constructionofW_yandW_qstsigma(W_y)=W_q}. According to $(\ref{eq:sigma^-1(p')capU_y=sigma_j^-1(p')j=1n_pq})$ we have for any $p'\in B_p$:
\begin{equation}\label{eq:expressionofv(q')onV_q}
    v_q(p')= \sum_{j=1}^{n_{p,q}}(h\circ\sigma_{j,p,q}^{-1})(p')d_{\sigma_{j,p,q}^{-1}(p')}\sigma\cdot(1,0).
\end{equation}
This formula (\ref{eq:expressionofv(q')onV_q}) shows that $v_q$ is holomorphic on $B_p$. Since these arguments are valid for any $p\in W_q\backslash E_R$, the vector field $v_q:W_q\backslash E_R\longrightarrow T\mathbb{P}^2$ is holomorphic.\\

\noindent\textit{2.} Let us fix $p\in W_q\backslash E_R$ and let $\Phi:B_p\to (\cmplex^2)^{n_{p,q}}$ be defined by $\Phi:p'\in B_p\mapsto(\sigma_{j,p,q}^{-1}(p'))_{1\leq j\leq n_{p,q}}$. Let us defined $Z_{p,q}:=\{(Z_j,W_j)_j\in(\cmplex^2)^{n_{p,q}},\ \exists i<j : Z_i=Z_j\},$ there exists $p'\in B_p\backslash \Phi^{-1}(Z_{p,q})$. Let us write $\Phi(p')=(x_j')_j=(Z_j,W_j)_j$
and let $\vec V_j:=d_{x_j'}\sigma\cdot(1,0)$. Taking $\lambda_1:=1$ and $\lambda_j:=0,\ j\geq2$, we have $\sum_{j=1}^{n_{p,q}}\lambda_j\vec V_{j} = \vec V_1 \neq \vec 0$, since $x_1'$ is not a critical point of $\sigma$. Since $\Phi(p')\not\in Z_{p,q}$, there exists a polynomial function $P(Z)$ such that $P(Z_j)=\lambda_j$. Let $h(Z,W):=P(Z)$ on $\cmplex^2$, we deduce from (\ref{eq:expressionofv(q')onV_q}) that $v_q(p')=\sum_{j=1}^{n_{p,q}}P(Z_j)\times \vec V_j=\vec V_1\neq0$.\qed\\

We prove now that $v_q$ extends holomorphically through $E_R$. We recall that $E_R$ is included in an analytic subset of codimension $1$ of $\sigma(D_R)$:
\begin{equation}\label{eq:E_R'=f^m(Crit f^m)capsigma(D_R)}
    E_R\subset E_R'\ ,\ \mathrm{with}\ E'_R:=f^{n_R}(\mathrm{Crit}\ f^{n_R})\cap \sigma(D_R).
\end{equation}

\begin{prop} The vector field $v_q$ admits a unique holomorphic extension on $W_q$, still denoted $v_q$.
\end{prop}
\noindent\textbf{\underline{Proof :}} In this proof we extends $v_q$ through the analytic subset $E_R'$ defined by (\ref{eq:E_R'=f^m(Crit f^m)capsigma(D_R)}), it extends in particular $v_q$ through $E_R$.

The content of the proposition is of local nature, thus we can assume that there exists $\xi:W_q\to \xi(W_q)\subset\cmplex^2$ a chart, holomorphic on a neighborhood of $\overline{W}_q$.
We define $v_*:=\xi_*v_q$ on $\xi(W_q\backslash E_R')=\xi(W_q)\backslash F'$, with $F':=\xi(E_R'\cap W_q)$. Let us write $v_*=(\alpha,\beta)$.

Let us fix $p\in W_q\backslash E_R'$ and let us consider the ball $B_p\subset W_q\backslash E_R$ we have constructed in Lemma \ref{lemma:constructionofW_yandW_qstsigma(W_y)=W_q}. Up to reduce this ball we have $B_p\subset W_q\backslash E_R'$. For $p'\in B_p$, $v_*(\xi(p')) = \sum_{j=1}^{n_{p,q}}(h\circ \sigma_{j,p,q}^{-1}(p'))\times \left(d_{p'}\xi\circ d_{\sigma_{j,p,q}^{-1}(p')}\sigma\right)\cdot(1,0)$ by \eqref{eq:expressionofv(q')onV_q}, and $n_{p,q}\leq d^{2n_R}$. Thus we deduce that $||v_*||\leq d^{2n_R} M_{q}$ on $\xi(B_p)\backslash F'$, where $M_{q}:=\sup_{D_R}|h|\times\sup_{W_q}||d\xi||\times\sup_{D_R}||d\sigma\cdot(1,0)||$. 

$M_q$ is finite and independent of $p$, thus by Riemann extension theorem (cf. {\cite[Proposition 1.1.7]{huybrechts2005complex}}), $\alpha$ and $\beta$ admit unique extensions $\widetilde{\alpha}$ and $\widetilde{\beta}$ on $\xi(W_q)$. Finally, the vector field $\xi^*(\widetilde{\alpha},\widetilde{\beta})$ extends uniquely $v_q$ on $W_q$.\qed\\

The interest of this vector field $v_q$ is that for each $p\in\mathrm{Reg}_{\sigma_R}(J)\cap W_q$, it is tangent to the distribution of directions $\mathcal{D}_p$ on a neighborhood of $p$ included in $V_p\cap W_q$:

\begin{lemme}\label{lemme:v_q(t)inD(t)} Let $p\in\mathrm{Reg}_{\sigma_R}(J)\cap W_q$. Let $B_p$ be the ball constructed in $\mathrm{Lemma\ \ref{lemma:constructionofW_yandW_qstsigma(W_y)=W_q}}$. There exists a ball $B_{p,q}\subset B_p\cap V_p\cap W_q$ centered at $p$ such that:
$$\forall t\in B_{p,q},\ v_q(t)\in \mathcal{D}_p(t).$$
\end{lemme}
\noindent\textbf{\underline{Proof :}} According to (\ref{eq:expressionofv(q')onV_q}), $v_q(p')= \sum_{j=1}^{n_{p,q}}(h\circ\sigma_{j,p,q}^{-1})(p')d_{\sigma_{j,p,q}^{-1}(p')}\sigma\cdot(1,0)$ for any $p'\in B_p$.
Observe that, because $p\in\mathrm{Reg}_{\sigma_R}(J)$, the inverse branches $\sigma_{j,p,q}^{-1}$ which appear are some of those we have constructed in Section \ref{sec:inversebranchesaboveregularvaluesofsigmainsidethejuliaset} for the regular value $p$ of $\sigma$. So for any $j\in\{1,\cdots,n_{p,q}\}$ and any $p'\in B_p$
we have $\sigma_{j,p,q}^{-1}(p')\in\mathrm{IB}_{R}(p)$, 
we refer to Definition \ref{defn:InverseBranchesIB_R(p)}. Hence by definition of $\mathcal{D}_p$ in (\ref{eq:Directionoffoliationmap}), we have for any $p'\in B_p\cap V_p$ and for any $j\in\{1,\cdots,n_{p,q}\}$, $\left[d_{\sigma_{j,p,q}^{-1}(p')}\sigma\cdot(1,0)\right] = \mathcal{D}_p(p')$ and thus $v_q(p')\in\mathcal{D}_p(p')$. If $B_{p,q}$ is a ball centered at $p$ and included in $B_p\cap V_p\cap W_q$, the conclusion follows.\qed\\

Up to a reduction of the neighborhood $U_{y,q}$ of $y$ in Lemma \ref{lemma:constructionofW_yandW_qstsigma(W_y)=W_q}, we can assume that $W_q=\sigma(U_{y,q})$ is equipped with holomorphic coordinates. As explained in Section \ref{sec:holomorphicfoliations}, we can use these local coordinates to create a non trivial holomorphic $1-$form $\omega_q$ on $W_q$ satisfying $\omega_q(v_q)\equiv0$.

\begin{defn}\label{defn:definitionofF_q} The (possibly singular) foliation induced by the $1-$form $\omega_q$ on $W_q$ is denoted $\mathcal{F}_q$.
\end{defn}

Observe that the equation $\omega_q(v_q)=0$ ensures that the leaves of the foliation $\mathcal{F}_q$ are tangent to the vector field $v_q$, at least where $v_q\neq0$. We refer to Section \ref{sec:holomorphicfoliations} for more details on foliations and vector fields. 

\begin{prop}\label{prop:patchingofFregwithF_qandF_qwithF_q'} \ 
\begin{enumerate}
    \item For any $q\in J\cap E_R$ and for any $p\in\mathrm{Reg}_{\sigma_R}(J)$, the foliations $\mathcal{F}_p$ and $\mathcal{F}_q$ coincide on each connected component $\mathcal{C}\subset V_p\cap W_q$ such that $\sigma_T(J\cap \mathcal{C})>0$.
    \item For any $(q,q')\in (J\cap E_R)^2$, the foliations $\mathcal{F}_q$ and $\mathcal{F}_{q'}$ coincide on each connected component $\mathcal{C}\subset W_q\cap W_{q'}$ such that $\sigma_T(J\cap \mathcal{C})>0$.
    \item For $q\in J\cap E_R$, $\sigma^*\mathcal{F}_q=\mathcal{F}_w$ on $U_{y,q}$, where $U_{y,q}$ is defined in $\mathrm{Lemma\ \ref{lemma:constructionofW_yandW_qstsigma(W_y)=W_q}}$.
\end{enumerate}
\end{prop}
\noindent\textbf{\underline{Proof :}} \ \\
\noindent\textit{1.} 
Since $\sigma_T(J\cap\mathcal{C})>0$ and $\sigma_T(E_R)=0$ (recall $\sigma_T(E)=0$, cf. Lemma \ref{prop:thedevolepingmap}), there exists $\widetilde{p}\in(J\cap \mathcal{C})\backslash E_R$. Let $\widetilde{\mathcal{C}}$ be the connected component of $\mathcal{C}\cap V_p\cap V_{\widetilde{p}}$ containing $\widetilde{p}$. Since $\widetilde{p}\in\mathrm{Supp}(T)$ we have $\sigma_T(\widetilde{\mathcal{C}})>0$. So according to Proposition \ref{prop:generalcaseoffiberconnections}, we have $\mathcal{D}_p(t)=\mathcal{D}_{\widetilde{p}}(t)$ for $t\in\widetilde{\mathcal{C}}$. Moreover for $t\in B_{\widetilde{p},q}\backslash\{v_q=0\}$, $[v_q(t)]=\mathcal{D}_{\widetilde{p}}(t)$ by Lemma \ref{lemme:v_q(t)inD(t)}. We deduce that $T_t\mathcal{F}_q=[v_q(t)]=\mathcal{D}_{\widetilde{p}}(t)=\mathcal{D}_p(t)$ for $t\in\widetilde{\mathcal{C}}\cap B_{\widetilde{p},q}\backslash\{v_q=0\}$. But $\mathcal{D}_p(t)=T_t\mathcal{F}_p$ for $t\in V_p$, thus we have $T_t\mathcal{F}_q=T_t\mathcal{F}_p$ for all $t\in\widetilde{\mathcal{C}}\cap B_{\widetilde{p},q}\backslash\{v_q=0\}\subset\mathcal{C}$. By using the second item and the third item of Lemma \ref{lemma:F_1=F_2onD^2(1)thenF_1=F_2onD^2(2)}, we conclude that finally $\mathcal{F}_q=\mathcal{F}_p$ on $\mathcal{C}$.\\

\noindent\textit{2.} Let $\mathcal{C}$ be a connected component of $W_q\cap W_{q'}$ and let us assume that $\sigma_T(J\cap \mathcal{C})>0$. As in the preceding point, there exists $p\in(\mathrm{Supp}(\sigma_T|_J)\cap \mathcal{C})\backslash E_R$. The connected component $\mathcal{C}'$ of $\mathcal{C}\cap V_{p}$ containing $p$ satisfies $\sigma_T(J\cap \mathcal{C}')>0$ since $p\in\mathrm{Supp}(\sigma_T|_J)$ . According to the previous item (applied first to $q$ and $p$ and then secondly to $q'$ and $p$) we have $\mathcal{F}_q=\mathcal{F}_p=\mathcal{F}_{q'}$ on $\mathcal{C}'$. Using again Lemma \ref{lemma:F_1=F_2onD^2(1)thenF_1=F_2onD^2(2)} we have at last $\mathcal{F}_q=\mathcal{F}_{q'}$ on $\mathcal{C}$.\\

\noindent\textit{3.} We have $\sigma_T(J\cap W_q)>0$ since $q\in J=\mathrm{Supp}(\sigma_T|_J)$. So there exists $p\in\mathrm{Reg}_{\sigma_R}(J)\cap W_q$. Let $B_{p,q}$ be the ball constructed in Lemma \ref{lemme:v_q(t)inD(t)}. Let $\sigma_{j,p,q}^{-1}$ be an inverse branch of $\sigma$ which appears in the formula (\ref{eq:expressionofv(q')onV_q}) and such that $x_j:=\sigma_{j,p,q}^{-1}(p)$ belongs to $\mathrm{IB}_R(p)\cap U_{y,q}$. Observe that $\sigma_{j,p,q}^{-1}(p)\in \mathrm{IB}_R(p)$ by definition of $\mathrm{IB}_R(p)$. The fact that we can choose $j$ such that $x_j\in U_{y,q}$ comes from the fact that $\sigma(U_{y,q})=W_q$, see Lemma \ref{lemma:constructionofW_yandW_qstsigma(W_y)=W_q}. 

According to the first item, the foliations $\mathcal{F}_p$ and $\mathcal{F}_q$ coincide on the ball $B_{p,q}\subset V_p\cap W_q$ (since $\sigma_T(J\cap B_{p,q})>0$ because $p\in J=\mathrm{Supp}(\sigma_T|_J)$). Thus we have $\sigma^*\mathcal{F}_q=\sigma^*\mathcal{F}_p\ \mathrm{on}\ \sigma^{-1}(B_{p,q})$.
But we have also $\sigma^*\mathcal{F}_p=\mathcal{F}_w\ \mathrm{on}\ \mathrm{IB}_R(p)$, see Definition \ref{defn:definitionofF_pbyomega_pandsigma^*F=F_wonIB_R(p)}. So we deduce that $\sigma^*\mathcal{F}_q=\mathcal{F}_w\ \mathrm{on}\ \sigma^{-1}(B_{p,q})\cap\mathrm{IB}_R(p),$
and the point $x_j$ belongs to $\sigma^{-1}(B_{p,q})\cap\mathrm{IB}_R(p)\cap U_{y,q}$. 
So $\sigma^*\mathcal{F}_q=\mathcal{F}_w$ on $U_{y,q}$ by analytic continuation (see Lemma \ref{lemma:F_1=F_2onD^2(1)thenF_1=F_2onD^2(2)}). \qed

\subsubsection{Step 4: Finite covering of \texorpdfstring{$J$}{TEXT} by foliated neighborhoods}\label{sec:finitecoveringofJ}

We have constructed for each point $p\in J$ a neighborhood $\mathcal{O}_p$ of $p$ ($\mathcal{O}_p:=V_p$ if $p\in\mathrm{Reg}_{\sigma_R}(J)$ or $\mathcal{O}_p=W_p$ if $p\in J\cap E_R$), equipped with a foliation $\mathcal{F}_p$ on $\mathcal{O}_p$. Propositions \ref{prop:T_p'P^2=D(p')} and \ref{prop:patchingofFregwithF_qandF_qwithF_q'} ensure that two foliations $\mathcal{F}_p$ and $\mathcal{F}_q$ of this collection coincide on $\mathcal{O}_p\cap \mathcal{O}_q$ when the connected components of $\mathcal{O}_p\cap\mathcal{O}_q$ are charged by $\sigma_T$ or $\sigma_T|_J$. Let us modify the covering $J\subset\cup_{p\in J}\mathcal{O}_p$ such that the two following properties hold: when two open sets of the covering intersect themselves, the foliations coincide on the intersection ; and every connected component of the covering contains an element of $J$. First, by an argument of Riemannian geometry, if we reduce the neighbourhoods $\mathcal{O}_p$ sufficiently, they are connected and they satisfy the following property: if two sets $\mathcal{O}_p$ and $\mathcal{O}_q$ intersect themselves then $\mathcal{O}_p\cap\mathcal{O}_q$ is also connected. Second, extract by compactness of $J$ a finite covering $J\subset\cup_{k=1}^{L}\mathcal{O}'_{p_k}$, where $\mathcal{O}'_{p_k}$ is a neighborhood of $p_k$ such that $\overline{\mathcal{O}'_{p_k}}\subset\mathcal{O}_{p_k}$. Third, for each $p_k$ define $\widetilde{\mathcal{O}_{p_k}}$ the open subset of $\mathcal{O}'_{p_k}$ obtained by depriving the sets $\overline{\mathcal{O}'_{p_k}\cap\mathcal{O}'_{p_l}},\ 1\leq l\leq L,$ which do not intersect $J$. Finally, we assert that the following covering has the two desire properties:
\begin{equation}\label{eq:coveringofJbyfoliatedballs}
    J\subset \mathcal{V}:=\widetilde{\mathcal{V}}\ \backslash\left\{\mathrm{connected\ components\ that\ do\ not\ contain\ elements\ of}\ J\right\}
\end{equation}
where $\widetilde{\mathcal{V}}:={\underset{k=1}{\overset{L}{\cup}}}\ \widetilde{\mathcal{O}_{p_k}}=\left({\underset{i=1}{\overset{N}{\cup}}}\ \widetilde{V}_{p_i}\right)\cup\left({\underset{j=1}{\overset{M}{\cup}}}\ \widetilde{W}_{q_j}\right).$ Indeed, if $\widetilde{\mathcal{O}}_{p_k}\cap\widetilde{\mathcal{O}}_{p_l}\neq\emptyset$ it means that $\overline{\mathcal{O}'_{p_k}\cap\mathcal{O}'_{p_l}}\cap J\neq\emptyset$, and thus $\sigma_T(J\cap \mathcal{O}_{p_k}\cap\mathcal{O}_{p_l})>0$ since $J=\mathrm{Supp}(\sigma_T|_{J})$. Using Proposition \ref{prop:generalcaseoffiberconnections} or \ref{prop:patchingofFregwithF_qandF_qwithF_q'} we conclude that $\mathcal{F}_{p_k}=\mathcal{F}_{p_l}$ on $\mathcal{O}_{p_k}\cap\mathcal{O}_{p_l}$, since this intersection is connected.

In particular, the foliations $(\mathcal{F}_{p_i})_{i\in\{1,\cdots,N\}}$ and $(\mathcal{F}_{q_j})_{j\in\{1,\cdots,M\}}$ patch all together on $\widetilde{\mathcal{V}}$ and form a (possibly singular) foliation:

\begin{defn}\label{defn:thefoliationFextractedfromthecoveringofJ} We denote 
\begin{equation*}
    \mathcal{F}:=\left(\bigcup_{i=1}^N\mathcal{F}_{p_i}|_{\widetilde{V}_{p_i}}\right)\cup \left(\bigcup_{j=1}^M\mathcal{F}_{q_j}|_{\widetilde{W}_{q_j}}\right)
\end{equation*} 
the foliation on the covering $\widetilde{\mathcal{V}}$ given by this patching.
\end{defn}

To finish the proof of Theorem \ref{thm:foliatedJ'} it remains to show that $\sigma^*\mathcal{F}=\mathcal{F}_w$ on $\sigma^{-1}(\mathcal{V})$.
It is done in the next section and use the Proposition \ref{rmk:tautologicalremark} below. Let us introduce an open set $\mathcal{W}\subset\mathbb{C}^2$ such that $\sigma(\mathcal{W})=\widetilde{\mathcal{V}}$ and on which we are able to compute $\sigma^*\mathcal{F}$. We denote for each $p_i$ and $q_j$:
$$\widetilde{\mathrm{IB}}_{R}(p_i) := \mathrm{IB}_R(p_i)\cap \sigma^{-1}\left(\widetilde{V}_{p_i}\right)\ \mathrm{and}\ \widetilde{U}_{y_j,q_j}:= U_{y_j,q_j}\cap \sigma^{-1}\left(\widetilde{W}_{q_j}\right),$$
where $U_{y_j,q_j}$ is defined in Lemma \ref{lemma:constructionofW_yandW_qstsigma(W_y)=W_q}.
The open set $\mathcal{W}$ is then given by:
\begin{equation*}
    \mathcal{W} := \left(\bigcup_{i=1}^{N} \widetilde{\mathrm{IB}}_{R}(p_i)\right)\cup \left(\bigcup_{j=1}^{M}\widetilde{U}_{y_j,q_j}\right).\\
\end{equation*}
For each $i$ and $j$ one has $\sigma\left(\widetilde{\mathrm{IB}}_R(p_i)\right)=\widetilde{V}_{p_i}$ and $\sigma\left(\widetilde{U}_{y_j,q_j}\right)=\widetilde{W}_{q_j}$, and thus
\begin{equation}\label{eq:sigma(W)=V}
    \sigma(\mathcal{W})=\widetilde{\mathcal{V}}.
\end{equation}
By using Definition \ref{defn:definitionofF_pbyomega_pandsigma^*F=F_wonIB_R(p)} and of the third item of Proposition \ref{prop:patchingofFregwithF_qandF_qwithF_q'}, we observe that:
\begin{prop}\label{rmk:tautologicalremark} We have $\sigma^*\mathcal{F}=\mathcal{F}_w\ \mathrm{on}\ \mathcal{W}.$
\end{prop}

\subsubsection{Step 5 : The foliation is horizontal on \texorpdfstring{$\sigma^{-1}(\mathcal{V})$}{TEXT}}\label{sec:invarianceofthefoliationF}

This is the last step to complete the proof of Theorem \ref{thm:foliatedJ'}.

\begin{prop}\label{lemma:liftofcontinuouspathbysigma} Let $\mathcal{C}$ be a connected component of $\sigma^{-1}(\mathcal{V})$. Then there exists $x\in \mathcal{C}\backslash(\mathrm{Crit}\ \sigma)$ such that $\sigma(x)\in J\backslash E_R$.
\end{prop}
\noindent\textbf{\underline{Proof :}} There exists $R'\geq R$ large enough such that $\mathcal{C}\cap \mathbb{D}^2_{R'}\neq\emptyset$. Assume $R'=R$ for simplicity. The set $\Sigma_R:=D_R\cap \sigma^{-1}(E_R)$ is an analytic subset of $D_R$ of codimension $\geq1$, thus $D_R\backslash\Sigma_R$ is connected and $(\mathcal{C}\cap D_R)\backslash\Sigma_R\neq\emptyset$. So if $x_0\in (\mathcal{C}\cap D_R)\backslash\Sigma_R$ is an arbitrary point, there exists $U_R\subset D_R\backslash\Sigma_R$ an open connected set, relatively compact in $D_R\backslash\Sigma_R$, which contains $x_0$. The point $\sigma(x_0)$ belongs to $\mathcal{V}$, let $\mathcal{V}_0\subset\mathcal{V}$ be the connected component containing $\sigma(x_0)$. 

By construction \eqref{eq:coveringofJbyfoliatedballs} of $\mathcal{V}$, we have $\mathcal{V}_0\cap J\neq\emptyset$. Since $J=\mathrm{Supp}(\sigma_T|_J)$ and since $\sigma_T(E_R)=0$, 
we then deduce the existence of a point $p\in \mathcal{V}_0\cap J\backslash E_R$. 
We can therefore increase $U_R$ if necessary to ensure that $\sigma(U_R)\cap \mathcal{V}_0$ contains a Lipschitz path $\gamma:[0,1]\to\sigma(U_R)\cap\mathcal{V}_0$ such that $\gamma(0)=\sigma(x_0)$ and $\gamma(1)=p$. Then according to Corollary \ref{cor:relevementdecheminparsigma}, there exists $\widetilde{\gamma}:[0,1]\to U_R$ a continuous path such that $\sigma\circ\Tilde{\gamma}=\gamma$ and $\Tilde{\gamma}(0)=x_0$. Since $\widetilde{\gamma}([0,1])\subset \sigma^{-1}(\mathcal{V})$ with $\widetilde{\gamma}(0)\in\mathcal{C}$, and since $\mathcal{C}$ is a connected component of $\sigma^{-1}(\mathcal{V})$, we must have $\widetilde{\gamma}([0,1])\subset \mathcal{C}$. In particular, we have $x:=\Tilde{\gamma}(1)\in \mathcal{C}\backslash(\mathrm{Crit}\ \sigma)$ such that $\sigma(x)=\gamma(1)=p\in J\backslash E_R$.\qed

\begin{lemme}\label{lemma:sigma^*F=F_w} The foliation $\sigma^*\mathcal{F}|_{\sigma^{-1}(\mathcal{V})}$ coincide with the horizontal foliation $\mathcal{F}_w$. 
\end{lemme}
\noindent\textbf{\underline{Proof :}} Let $\mathcal{C}$ be a connected component of $\sigma^{-1}(\mathcal{V})$. According to Proposition \ref{lemma:liftofcontinuouspathbysigma} there exists $x\in\mathcal{C}\backslash(\mathrm{Crit}\ \sigma)$ such that $p:=\sigma(x)\in J\backslash E_R$. Thus there exists $U_{x}\subset\mathcal{C}$ a connected open neighborhood of $x$ such that $\sigma_{x}:=\sigma|_{U_{x}}:U_{x}\to \sigma(U_{x})=:V_{x}$ is a biholomorphism. By (\ref{eq:sigma(W)=V}) there exists $y\in \mathcal{W}$ such that $\sigma(y)=p$. Since $y\in D_R$ and since $p$ is not a critical value of $\sigma|_{D_R}$, $y$ is not a critical value of $\sigma$. So let $U_y\subset \mathcal{W}$ be a open neighborhood of $y$ such that $\sigma_y:=\sigma|_{U_y}:U_y\to\sigma(U_y)=:V_y$ is a biholomorphism. 

We can assume that $V_{x}\subset V_y$ and $\phi:=\sigma_{y}^{-1}\circ\sigma_{x}:U_{x}\longrightarrow \sigma_{y}^{-1}(V_{x})$ is well defined. Observe that $\sigma_T(V_{x})>0$ (since $p\in J$) and that $V_{x}$ is connected. So according to the first item of Proposition \ref{prop:generalcaseoffiberconnections}, the map $\phi$ has the form $\phi(z,w)=(A(z,w),B(w))$ and it preserves $\mathcal{F}_w$. Observe that $\phi^*(\sigma_y^*\mathcal{F}|_{V_{x}})=\sigma_{x}^*\mathcal{F}|_{V_{x}}$ by definition of $\phi$. Moreover, according to Proposition \ref{rmk:tautologicalremark}, we have $\sigma_y^*(\mathcal{F}|_{V_y})=\mathcal{F}_w|_{U_y}$, thus we have $(\sigma_{x}^*\mathcal{F})|_{U_{x}}=(\phi^*\mathcal{F}_w)|_{U_{x}}=\mathcal{F}_w|_{U_{x}}$. We deduce that $\sigma^*\mathcal{F}=\mathcal{F}_w$ on $\mathcal{C}$ by analytic continuation (Lemma \ref{lemma:F_1=F_2onD^2(1)thenF_1=F_2onD^2(2)}). Since it is true for any connected component $\mathcal{C}$ the conclusion follows.\qed

\begin{otherlanguage}{english}
\bibliographystyle{abbrv}
\bibliography{biblio} 
\end{otherlanguage}

$ $ \\
\noindent {\footnotesize V. Tapiero}\\
{\footnotesize Universit\'e de Rennes}\\
{\footnotesize CNRS, IRMAR - UMR 6625}\\
{\footnotesize F-35000 Rennes, France}\\
{\footnotesize virgile.tapiero@univ-rennes.fr}\\

\end{document}